\documentclass[a4paper,12pt]{article}
\usepackage{amsmath,amssymb,amsthm,latexsym,graphicx}
\usepackage{marvosym}
\usepackage{dictsym}
\usepackage{natbib}
\usepackage{algorithmic}
\usepackage{longtable}

\theoremstyle{definition}

\newtheorem{lemma}{Lemma}
\newtheorem{proposition}[lemma]{Proposition}
\newtheorem{theorem}[lemma]{Theorem}

\newtheorem{definition}[lemma]{Definition}
\newtheorem{example}[lemma]{Example}
\newtheorem{remark}[lemma]{Remark}

\renewenvironment{proof}{\textbf{Proof:}}{}

\title{Lattice polygons and families of curves on rational surfaces}
\author{Niels Lubbes and Josef Schicho \thanks{This research was supported by the Austrian Science Fund (FWF)  in the frame of the  research project P-21461.}   }

\date{\today}

\begin{document} \maketitle \tableofcontents
\begin{abstract}
 First we solve the problem of finding minimal degree families on toric surfaces by reducing it to lattice geometry. Then we describe how to find minimal degree families on, more generally, rational complex projective surfaces. 
\end{abstract}

\section{Introduction}
 
\label{sec:intro_opt_fam} Every algebraic surface in projective space ${\textbf{P}}^r$ can be generated by a family of curves in projective space (e.g. the hyperplane sections). For a fixed surface, this can be done in infinitely many ways. Maybe the simplest family of algebraic curves is one where the curves have minimal genus, and among those one with minimal degree. In this paper, we study the families of genus zero curves of minimal degree, in the case where the given surface is rational (\textsection\ref{sec:fam_on_rat}). Classical examples of such families are the families of lines on a ruled surface -- with the single example of the nonsingular quadric in ${\textbf{P}}^3$ having two such families -- and the families of conics on a non-ruled conical surface (the surfaces with more than one family of conics have been classified in \cite{Schicho:00a}). \\[2mm]
 The paper, starts with a seemingly quite different topic, namely the study of discrete directions which minimize the width of a given convex lattice polytope (\textsection\ref{sec:clp}). As the lattice points reminds of sticks in a vineyard, we call the problem of finding all these directions the ``vineyard problem''; for the minimal directions, most sticks are aligned with others and one ``sees'' only a minimal number. We give an elementary solution, based on the notion of the adjoint lattice polytope, which is defined as the convex hull of the interior lattice points (see \textsection\ref{sec:va_on_vy}). \\[2mm]
 The vineyard problem is equivalent to the specialization of the problem of finding toric families of minimal degree on a given toric surface  (see Proposition~\ref{prop:tor_va_rel}). The main result of this paper is the fact that our elementary solution can be translated into the language of toric geometry, and then generalizes in a natural way so that it makes it possible to construct all minimal degree families of rational curves on arbitrary rational surfaces! In \textsection\ref{sec:fam_on_mprs}, we give a proof in the language of algebraic geometry (which subsumes then the elementary proof in \textsection\ref{sec:va_on_vy}). The methods are quite different, but, as the reader may check, there is a close analogy in the structure of the two proofs. \\[2mm]
 The algebraic geometry analogue of the adjoint lattice polytope is adjunction; this has been observed in \cite{ful1} (see also \cite{Schicho:03a, Schicho:07f}).

\subsection{Overview}
 The following table gives the problems and their solutions which are treated in this document:\\[2mm]
 {\tiny \begin{tabular}
{|l|l|l|} \hline problem & solution & description \\
\hline Definition~\ref{def:va_problem} & Theorem~\ref{thm:va_opt} & vineyard problem (or viewangle problem on vineyards) \\
\hline Definition~\ref{def:toric_problem} & Proposition~\ref{prop:tor_va_rel} & toric family problem on toric surfaces \\
\hline Definition~\ref{def:fam_mprs_problem} & Theorem~\ref{thm:optfam} & rational family problem on polarized rational surfaces \\
\hline Definition~\ref{def:fam_rat_problem} & Proposition~\ref{prop:optfam_rat} & rational family problem on rational surfaces \\
\hline 
\end{tabular}
 } \\[2mm]
 The second problem is reduced to the first problem and the fourth problem is  reduced to the third problem. In section \textsection\ref{sec:clp} we define convex lattice polygons and their adjoints. In \textsection\ref{sec:fam} we will define what we mean by family and give properties, of which Proposition~\ref{prop:fam_rat_prop} is most important. For \textsection\ref{sec:fam_on_tor} only Definition~\ref{def:fam} is needed of \textsection\ref{sec:fam}. In \textsection\ref{sec:adjoint_chain} we summarize the notions of  minimally polarized rational surface (mprs for short), adjoint relation  and adjoint chain, which are used in \textsection\ref{sec:fam_on_mprs}. See Remark~\ref{rem:analogy} for the analogy between \textsection\ref{sec:va_on_vy} and \textsection\ref{sec:fam_on_mprs}. 

\subsection{Guide for reading}
 
We explain the structure of this document.
The \textit{main-claims} are labeled by `\textbf{[a-z])}'.
A \textit{claim} is given by the sentence starting with `\textit{Claim [1-10]:}' and
is a step for proving the main claims.
The proof of a claim is given by the remaining sentences in the same paragraph.
We define each sentence in the proof of a claim to be a \textit{sub-claim}.

\section{Convex lattice polygons}
 
\label{sec:clp}
\begin{definition}
\label{def:lattice}
\textrm{\textbf{(lattice and dual lattice)}} 
  A \textit{lattice} $\Lambda_n$ is defined as ${\textbf{Z}}^n\subset{\textbf{R}}^n$.  Its \textit{dual lattice} $\Lambda_n^*$ is defined as $\textrm{Hom}_{{\textbf{Z}}}(\Lambda_n,{\textbf{Z}})$.  A \textit{lattice equivalence} is a map (translation, rotation, shearing and reflection):  \[\Phi: {\textbf{R}}^n \to {\textbf{R}}^n,\quad \overrightarrow{x} \mapsto A\overrightarrow{x} + \overrightarrow{y} \] where $A \in GL_n({\textbf{Z}})$ and $\overrightarrow{y} \in {\textbf{Z}}^n$.  We will denote $\Lambda_2$ by $\Lambda$. 
 
\end{definition}

\begin{definition}
\label{def:clp}
\textrm{\textbf{(convex lattice polygon)}} 
   Let $\Lambda$ be a two dimensional lattice. A \textit{convex lattice polygon} $\Gamma$ is the convex hull of a finite non-empty set of lattice points in $\Lambda$. Polygons are considered equivalent when they are lattice equivalent. 
 
\end{definition}

\begin{definition}
\label{def:clp_attr}
\textrm{\textbf{(attributes of polygons)}} 
  Let $\Gamma$ be a lattice polygon with lattice $\Lambda$. We call $\Gamma$ a \textit{shoe polygon}  if and only if  $\Gamma=\Box_{l,m,n}$ where  \[ \Box_{l,m,n}:=\textrm{ConvexHull}( (0,0),(0,l),(m,l),(m+n,0) ) \] where $l,m,n\in{\textbf{Z}}_{\ge 0}$ (see Figure~\ref{fig:clp_attr}.a)).  
\begin{figure}[h!]
 \centering \begin{tabular}
{|c|c|c|} {\includegraphics[width=6cm,height=3cm]{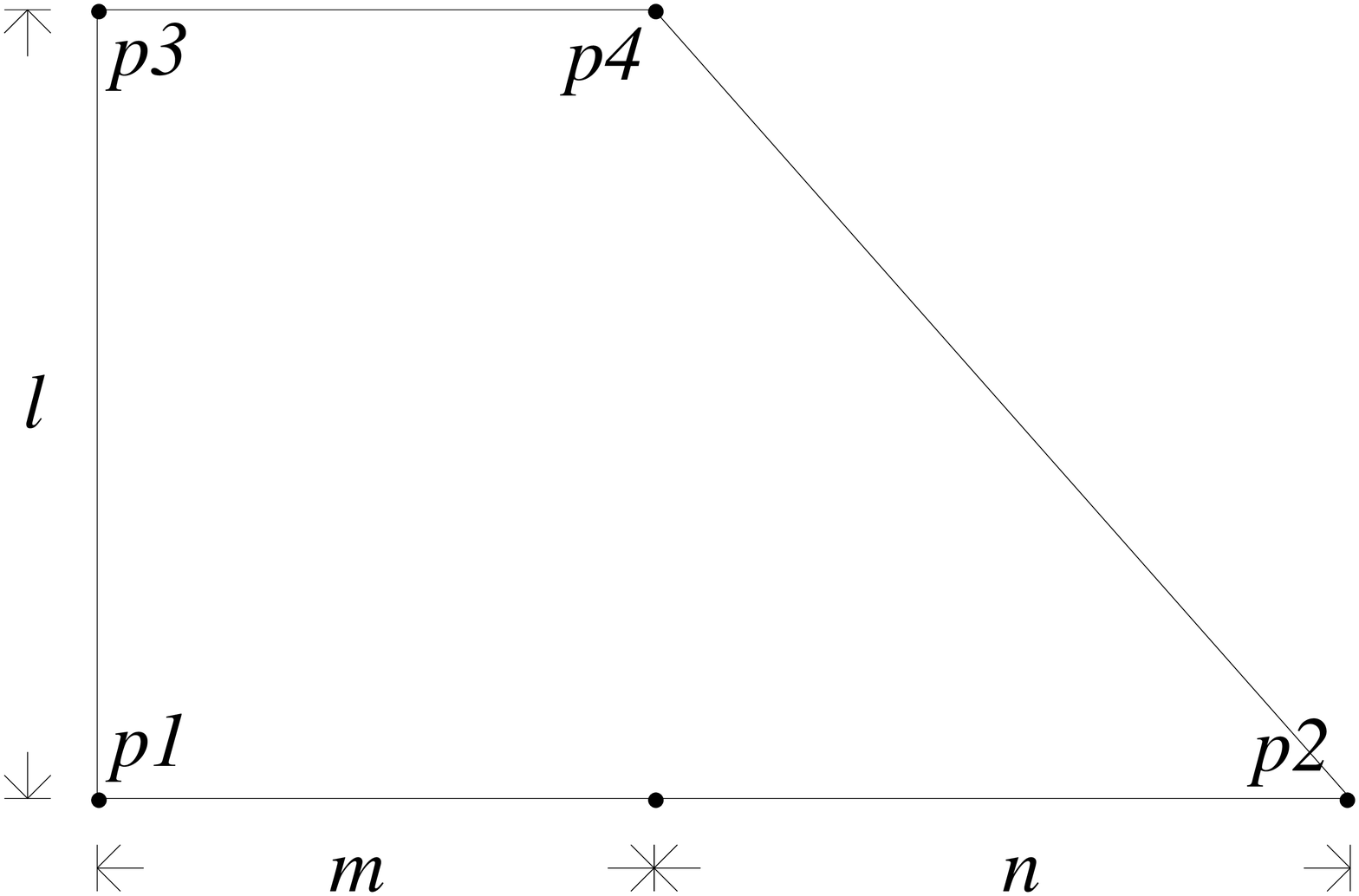}} & {\includegraphics[width=3cm,height=3cm]{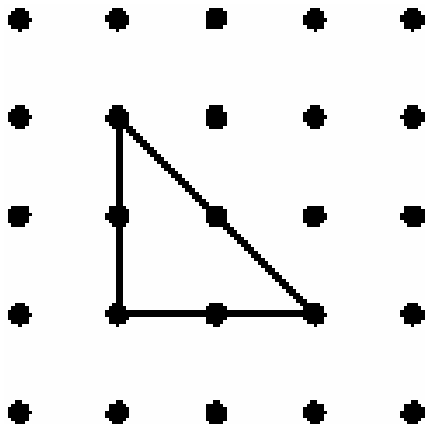}}   & {\includegraphics[width=3cm,height=3cm]{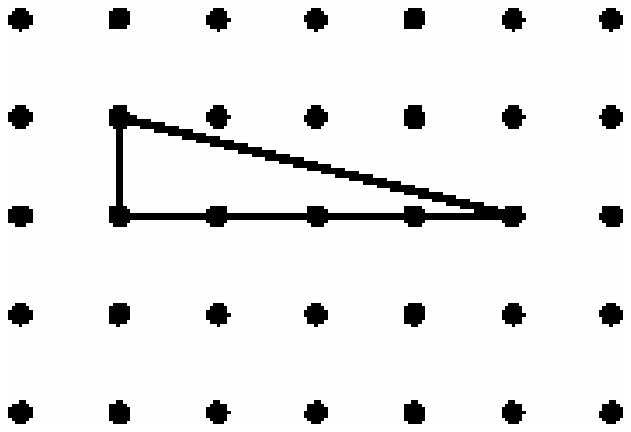}} \\
 a: $\Box_{l,m,n}$       & b: $\Box_{2,0,2}$             & c: $\Box_{1,0,>1}$ 
\end{tabular}
 \caption{convex lattice polytopes corresponding to $\Box_{l,m,n}$.} 
\label{fig:clp_attr} 
\end{figure} We call $\Gamma$ a \textit{standard triangle}  if and only if  $\Gamma=\Box_{l,0,l}$ with $l>0$ ( for example  Figure~\ref{fig:clp_attr}.b). We call $\Gamma$ a \textit{thin triangle}  if and only if  $\Gamma=\Box_{1,0,l}$ with $l>1$. ( for example  Figure~\ref{fig:clp_attr}.c). We call $\Gamma$ \textit{minimal}  if and only if    
 $\Gamma$ is not a point and either 
 $\Gamma$ has 1 interior lattice point or 
 $\Gamma$ has no interior lattice points. 

\end{definition}

\begin{definition}
\label{def:adjoint_polygon}
\textrm{\textbf{(adjoint polygon)}} 
   Let $\Gamma$ be a convex lattice polygon with lattice $\Lambda$. The \textit{adjoint polygon} $\Gamma'$ of $\Gamma$ is defined as the convex hull of the interior lattice points of $\Gamma$ (if there exist any). We denote the adjoint of $\Gamma$ taken $i$ times by $\Gamma^i$. 
 
\end{definition}

\begin{definition}
\label{def:va}
\textrm{\textbf{(viewangles and width)}} 
 Let $\Gamma$ be a convex lattice polygon with lattice $\Lambda$. A \textit{viewangle} for $\Gamma$ is a nonzero vector $h \in \Lambda^*-\{0\}$ in the dual lattice. The \textit{viewangle width} of a viewangle $h$ for $\Gamma$ is: \[ \textrm{width}_{\Gamma}: \Lambda^*\to{\textbf{Z}},~ h \mapsto \displaystyle\max_{v \in \Gamma}~h(v) - \displaystyle\min_{w \in \Gamma}~h(w). \] The \textit{width of a convex lattice polygon}  is the smallest possible viewangle width: \[ v(\Gamma)=\displaystyle\min_{h \in \Lambda^*-{0} }~\textrm{width}_{\Gamma}(h). \] The \textit{set of optimal viewangles} on $\Gamma$ is defined as  \[ S(\Gamma)=\{~ h \in \Lambda^*-\{0\} ~~|~~ \textrm{width}_\Gamma(h)=v(\Gamma) ~\}. \]
 
\end{definition}

\begin{definition}
\label{def:va_attr}
\textrm{\textbf{(attributes of viewangles)}} 
  Let $\Gamma$ be a convex lattice polygon with lattice $\Lambda$. Let $h \in \Lambda^*-\{0\}$ be a viewangle. \textit{\textbf{tight viewangles}}: 
  We call $h$ \textit{max-tight} for $\Gamma$  if and only if  $\Gamma'$ is defined and $\displaystyle\max_{v \in \Gamma}~h(v)=\displaystyle\max_{w \in \Gamma'}~h(w)+1$. We call $h$ \textit{min-tight} for $\Gamma$  if and only if  $\Gamma'$ is defined and $\displaystyle\min_{v \in \Gamma}~h(v)=\displaystyle\min_{w \in \Gamma'}~h(w)-1$. We call $h$ is \textit{tight} for $\Gamma$  if and only if  $h$ is max-tight and min-tight for $\Gamma$. 
    \textit{\textbf{edge viewangles}}: 
  We call $h$ a \textit{max-edge} for $\Gamma$  if and only if  $h(v)=h(w)=\displaystyle\max_{u \in \Gamma}~h(u)$ for some $v,w\in\Gamma$ where $v\neq w$. We call $h$ a \textit{min-edge} for $\Gamma$  if and only if  $h(v)=h(w)=\displaystyle\min_{u \in \Gamma}~h(u)$ for some $v,w\in\Gamma$ where $v\neq w$. We call $h$ an \textit{edge} for $\Gamma$  if and only if  $h$ is a max-edge and min-edge for $\Gamma$.

\end{definition}

\section{Minimal width viewangles for convex lattice polygons}
 
\label{sec:va_on_vy}
\begin{definition}
\label{def:va_problem}
\textrm{\textbf{(vineyard problem)}} 
  Given a convex lattice polygon $\Gamma$  find the width $v(\Gamma)$ and all optimal viewangles $S(\Gamma)$ (see Definition~\ref{def:va}). 
 
\end{definition}

\begin{example}\label{ex:va_problem}\textrm{\textbf{(vineyard problem)}}

 Let $\Gamma$ be the convex lattice polygon as in Figure~\ref{fig:va_problem} with viewangles $h_0=(1,-1)$ and $h_1=(1,0)$.  The origin is defined by the interior lattice point of $\Gamma$. 
\begin{figure}[h!]
 \centering {\includegraphics[width=3cm,height=3cm]{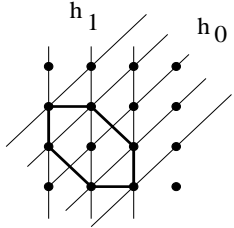}} \caption{A convex lattice polygon and two viewangles.} 
\label{fig:va_problem} 
\end{figure} 

 We have that $\textrm{width}_{\Gamma}(h_0)=4$ and $\textrm{width}_{\Gamma}(h_1)=2$. We find for this easy example that $v(\Gamma)=2$. The optimal viewangles are $h_1$, the horizontal viewangle $(0,1)$ and the diagonal viewangle $(-1,-1)$. 

\end{example}

\begin{lemma}\label{lem:va_low}\textrm{\textbf{(lowerbound)}}  Let $\Gamma$ be a convex lattice polygon which is not minimal  (see Definition~\ref{def:clp_attr} for minimal) with lattice $\Lambda$.  Let $h\in\Lambda^*-\{0\}$ be a viewangle. 
\begin{itemize} 
\item[] We have that $\textrm{width}_\Gamma(h)\ge \textrm{width}_{\Gamma'}(h)+2$, and equality holds  if and only if  $h$ is tight for $\Gamma$. 
\end{itemize}
\begin{proof}  
 Direct consequence of Definition~\ref{def:va} and Definition~\ref{def:va_attr}. 
 \end{proof} 
\end{lemma}

\begin{lemma}\label{lem:va_tight}\textrm{\textbf{(tight)}}  Let $\Gamma$ be a convex lattice polygon which is not minimal (see Definition~\ref{def:clp_attr} for minimal) with lattice $\Lambda$.  Let $h\in\Lambda^*-\{0\}$ be a viewangle. 
\begin{itemize} 
\item[\textbf{a)}] If $h$ is an edge of $\Gamma'$ then $h$ is tight for $\Gamma$. 
\item[\textbf{b)}] If $h$ is tight for $\Gamma'$ then $h$ is tight for $\Gamma$. 
\end{itemize}
\begin{proof} 
 We assume that $h$ is not max-tight for $\Gamma$ in the remainder of the proof.
 Let $p \in \Gamma$ be  such that  $\displaystyle\max_{v \in \Gamma }~h(v)=h(p)$.
 We will denote the lattice points in Figure~\ref{fig:tight} by the checkboard coordinates \textbf{a8} until \textbf{h1}. 
\begin{figure}[h!]
 \centering {\includegraphics[width=5cm,height=5cm]{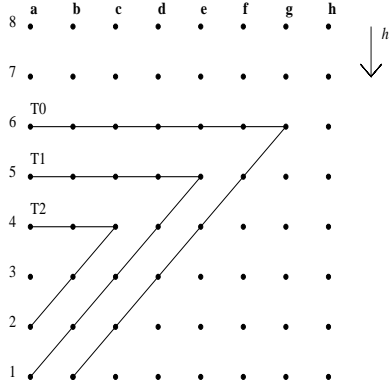}} \caption{Proof of Lemma~\ref{lem:va_tight}} 
\label{fig:tight} 
\end{figure}

\textit{Claim 1:} We may assume  without loss of generality  that $h=(1,0)$, $h(\textbf{e5})=\displaystyle\max_{v \in \Gamma'}~h(v)$ and $p$ is right of column \textbf{f}. 
  From the assumption that $h$ is not max-tight it follows that $p$ is right of column of \textbf{f}. 

 Let $S=\Gamma \cap L$ be a line segment where $L$ is the line corresponding to column \textbf{f}.

\textit{Claim 2:} The line segment $S$ doesn't contain interior lattice points and is not empty. 
  Suppose by contradiction that $S$ contains an interior lattice point $q$.  Then $q\in\Gamma'$ and $h(q)>h(\rm{\textbf{e5}})$. \textrm{\large\Lightning}

\textit{Claim 3:} We may assume  without loss of generality  that \textbf{f6} and \textbf{f5} are the lattice points above  respectively  under $S$. 
  From claim 2) it follows that $S$ is between \textbf{fi+1} and \textbf{fi} for some $i\in{\textbf{Z}}$.  We apply shearing such that \textbf{f6} and \textbf{f5} are the required points.  We have that $h$ remains unchanged under the corresponding dual transformation. 

 Let $Q=$ConvexHull(~\textbf{f6},~\textbf{f5}, $\Gamma~\cap$ the area right of column \textbf{f} ).

\textit{Claim 4:} The polygon $Q$ doesn't contain interior lattice points and is not empty. 
  It follows from the assumption that $h$ is not max-tight. 

 For example $Q$ is ConvexHull(\textbf{f6},\textbf{f5},\textbf{g7},\textbf{h7}) or ConvexHull(\textbf{f6},\textbf{f5},\textbf{h6}). For constructing examples it is required that $Q$ doesn't contain interior lattice points and that \textbf{e5} is between the line through (\textbf{f6},p) and the line through (\textbf{f5},p).
 Let $\widehat{\Gamma}=$ConvexHull(~$\Gamma - Q$,~\textbf{g6}~$)$.
 Let $T_0, T_1$ and $T_2$ be the area contained by the corresponding line as in Figure~\ref{fig:tight}.

\textit{Claim 5:} We have that $\widehat{\Gamma} \subseteq T_0,~ \Gamma'\subseteq T_1 {\textrm{~and~}}  \Gamma''\subseteq T_2$. 
  Suppose by contradiction that $\widehat{\Gamma}$ has a point outside of $T_0$.  It follows that $\Gamma$ is not convex. \textrm{\large\Lightning}   We have that $\Gamma' \subseteq \widehat{\Gamma}' \subseteq T_0'=T_1$ and $\Gamma''\subseteq T_1'=T_2$.

\textit{Claim 6:} If $h$ is not max-tight for $\Gamma$ then $h$ is not a max-edge of $\Gamma'$. 
  From  claim 5) and Figure~\ref{fig:tight} it follows that $h$ reaches the maximum only once for $\Gamma'\subseteq T_1$ at \textbf{e5}.  It follows that $h$ is not a max-edge for $\Gamma'$.

\textit{Claim 7:} If $h$ is not max-tight for $\Gamma$ then $h$ is not max-tight for $\Gamma'$. 
  From claim 5) and Figure~\ref{fig:tight} it follows that $h$ reaches a maximum for $\Gamma''\subseteq T_2$ on or left of column \textbf{c}.  It follows that $h$ is not max-tight for $\Gamma'$.

\textit{Claim 8:} From claim 6) and claim 7) it follows that \textbf{a)} and \textbf{b)}. 
  The proof of claim 6) and claim 7) for min-edge and min-tight is completely symmetric.  The statements are dual to \textbf{a)} and \textbf{b)}. 

 \end{proof}
\end{lemma}

\begin{proposition}\label{prop:va_min}\textbf{\textrm{(classification of optimal viewangles for minimal convex lattice polygons)}}
\begin{itemize} 
\item[\textbf{a)}] All the optimal viewangles on minimal convex lattice polygons are classified in Figure~\ref{fig:va_opt}.
\item[\textbf{b)}] If $\Gamma$ has a thin triangle ( id est  Figure~\ref{fig:va_opt}.20) as adjoint then $S(\Gamma)=S(\Gamma')$ and the optimal viewangle is tight for $\Gamma$.
\begin{figure}[h!]
 \centering \begin{tabular}
{llll} {\includegraphics[width=2cm,height=2cm]{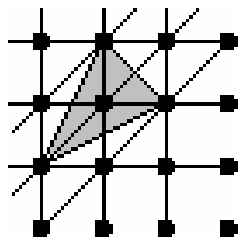}} & {\includegraphics[width=2cm,height=2cm]{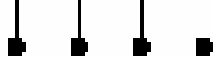}} & {\includegraphics[width=2cm,height=2cm]{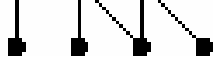}} & {\includegraphics[width=2cm,height=2cm]{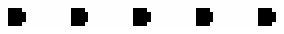}} \\
  1: $(2,3)$ tight &2: $(2,2)$ tight &3: $(2,3)$ tight &4: $(2,1)$ tight \\
 {\includegraphics[width=2cm,height=2cm]{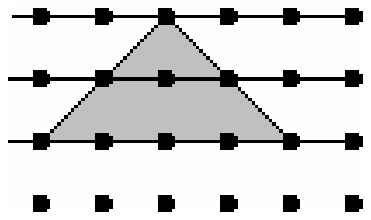}} & {\includegraphics[width=2cm,height=2cm]{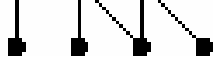}} & {\includegraphics[width=2cm,height=2cm]{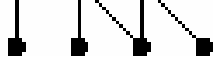}}  & {\includegraphics[width=2cm,height=2cm]{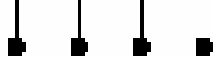}}  \\
   5: $(2,1)$ tight & 6: $(2,3)$ tight & 7: $(2,3)$ tight & 8: $(2,2)$ tight \\
 {\includegraphics[width=2cm,height=2cm]{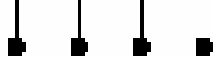}}  & {\includegraphics[width=2cm,height=2cm]{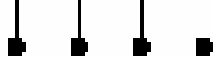}} & {\includegraphics[width=2cm,height=2cm]{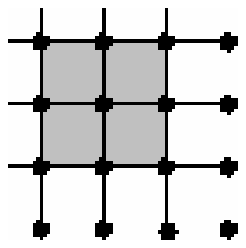}} & {\includegraphics[width=2cm,height=2cm]{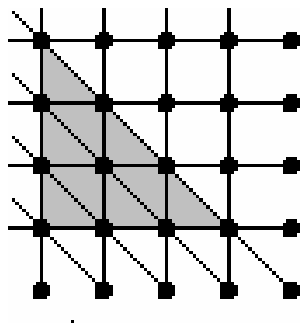}} \\
   9: $(2,2)$ tight &10: $(2,2)$ tight &11: $(2,2)$ tight &12: $(3,3)$ $\triangle_4$ \\
 {\includegraphics[width=2cm,height=2cm]{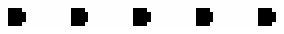}} & {\includegraphics[width=2cm,height=2cm]{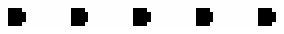}} & {\includegraphics[width=2cm,height=2cm]{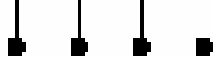}} & {\includegraphics[width=2cm,height=2cm]{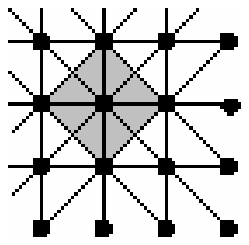}}  \\
  13: $(2,1)$  tight &14: $(2,1)$  tight &15: $(2,2)$  tight &16: $(2,4)$  tight \\
 {\includegraphics[width=2cm,height=2cm]{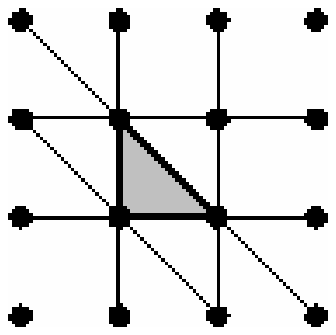}} & {\includegraphics[width=2cm,height=2cm]{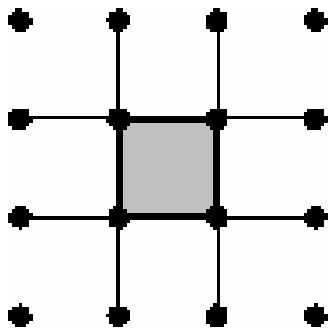}} & {\includegraphics[width=2cm,height=2cm]{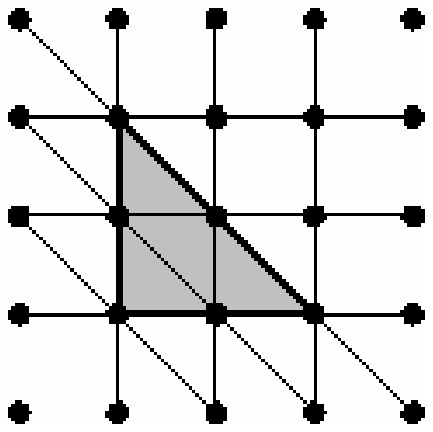}} & {\includegraphics[width=2cm,height=2cm]{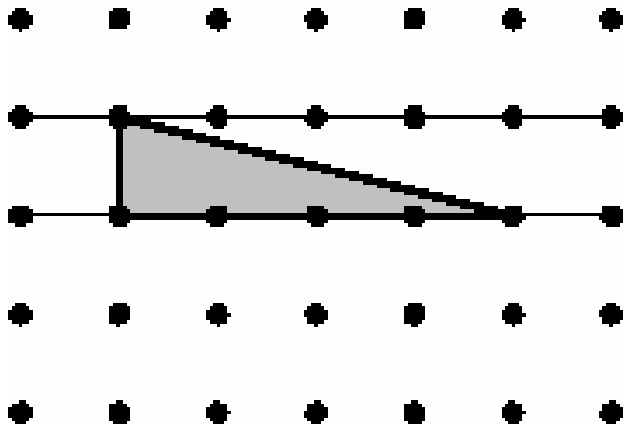}}  \\
  17: $(1,3)$  $\triangle_1$ &18: $(1,2)$  edge &19: $(2,3)$  $\triangle_2$ &20: $(1,1)$ \\
 {\includegraphics[width=2cm,height=2cm]{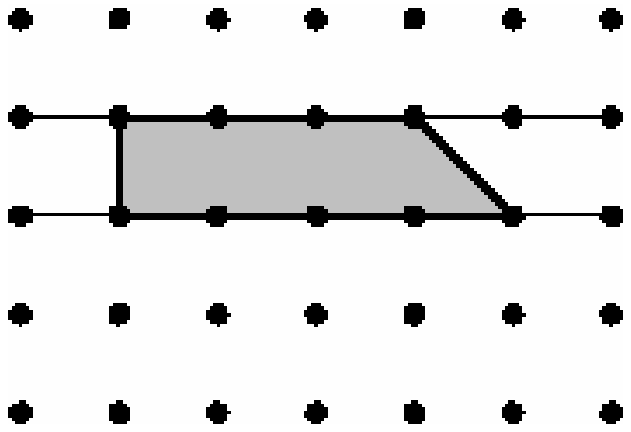}} & {\includegraphics[width=2cm,height=2cm]{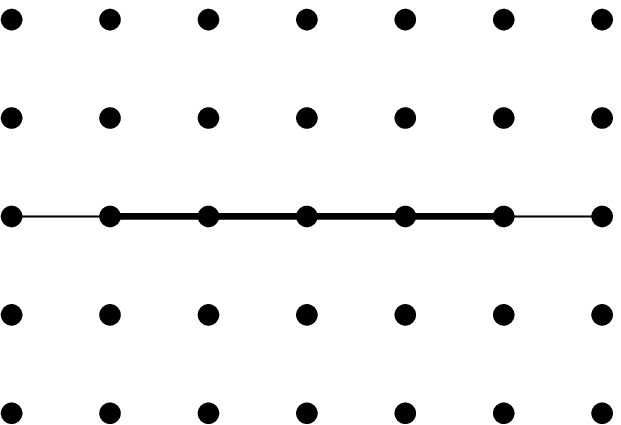}} & & \\
  21: $(1,1)$  edge &22: $(0,1)$  edge & & 
\end{tabular}
 \caption{ All the optimal viewangles for minimal convex lattice polygons and $(v(\Gamma),\#S)$ where $\#S$ is the number of optimal viewangles. We denote standard triangles of length $i$ by $\triangle_i$. } 
\label{fig:va_opt} 
\end{figure} 
\end{itemize}
\begin{proof} 
 The classification of minimal convex lattice polygons (see Definition~\ref{def:clp_attr}) can be found in \cite{Schicho:03a}.
 The classification of the optimal viewangles in Figure~\ref{fig:va_opt} is a direct result of tedious case by case inspection.
 Let's assume $\Gamma$ is a convex lattice polygon  such that  $\Gamma'=\Box_{1,0,l}$ and $l>1$ ( id est  thin triangle).

\textit{Claim:} We have $l=2$ and the optimal direction of $\Gamma$ is tight. 
  If $l>2$ then $\Gamma$ is not convex.  There are a finite number of possibities for $\Gamma$, each for which the optimal direction is tight. 

 \end{proof} 
\end{proposition}

\begin{definition}
\label{def:casedist}
\textrm{\textbf{(case distinction)}} 
   Let $\Gamma$ be a convex lattice polygon with lattice $\Lambda$. We distinguish between the following cases where $\Gamma$ is not minimal except at A0: \\[2mm]
 \begin{tabular}
{|l|l|l|} \hline    & $\Gamma$              & $\Gamma'$ \\
\hline A0 & minimal               & point or emptyset \\
\hline A1 & standard triangle     & standard triangle \\
\hline A2 & not standard triangle & standard triangle \\
\hline A3 & not standard triangle & minimal and not standard triangle \\
\hline A4 & not standard triangle & not minimal and not standard triangle \\
\hline 
\end{tabular}
 \\[2mm]
 See Definition~\ref{def:clp_attr} for the notion of standard triangle. 
 
\end{definition}

\begin{theorem}\label{thm:va_opt}\textbf{\textrm{(optimal viewangles)}}  Let $\Gamma$ be a convex lattice polygon with lattice $\Lambda$.  Let $S(\Gamma)$ be the set of all optimal viewangles of $\Gamma$.  Let A0 until A4 be as in Definition~\ref{def:casedist}. 
\begin{itemize} 
\item[\textbf{a)}] If A0 then  
 $S(\Gamma)$ and $v(\Gamma)$ are as classified in Figure~\ref{fig:va_opt}. 
\item[\textbf{b)}] If A1 then  
 $S(\Gamma)$ contains exactly its 3 edges and 
 $v(\Gamma)=v(\Gamma')+3$. 
\item[\textbf{c)}] If A2 then  
 $S(\Gamma)=\{~ h \in S(\Gamma') ~~|~~ h $ tight for $\Gamma ~\}$ and 
 $v(\Gamma)=v(\Gamma')+2$. 
\item[\textbf{d)}] If A3 or A4 then  
 $S(\Gamma)=S(\Gamma')$ and 
 $v(\Gamma)=v(\Gamma')+2$. 
\end{itemize}
\begin{proof}
 We have that \textbf{a)} and \textbf{b)} are a direct consequence of Proposition~\ref{prop:va_min} and the definition of the standard triangle.
 Let $T(\Gamma)=\{~ h ~~|~~ h\in S(\Gamma)$ and $h$ is tight $~\}$.

\textit{Claim 1:} If $T(\Gamma')\neq\emptyset$ then $S(\Gamma)=T(\Gamma)$. 
  From Lemma~\ref{lem:va_low} and Lemma~\ref{lem:va_tight}.a it follows that if $h\in T(\Gamma')$ then $\textrm{width}_\Gamma(h) = v(\Gamma')+2$.  From Lemma~\ref{lem:va_low} it follows that if $h\in S(\Gamma)$ then $\textrm{width}_\Gamma(h) \ge v(\Gamma')+2$ and equality holds  if and only if  $h \in T(\Gamma)$.

\textit{Claim 2:} If $S(\Gamma')=T(\Gamma')$ then $S(\Gamma)=S(\Gamma')$. 
  From Lemma~\ref{lem:va_low} and Lemma~\ref{lem:va_tight}.a it follows that if $h\in T(\Gamma')$ then $\textrm{width}_\Gamma(h) = v(\Gamma')+2$.  It follows that $T(\Gamma)\supseteq T(\Gamma')$.  If $h\in T(\Gamma)$ then $\textrm{width}_\Gamma(h) =\textrm{width}_{\Gamma'}(h)+2$ and thus $h\in S(\Gamma')$.  It follows that $T(\Gamma)\subseteq S(\Gamma')$.  From claim 2) and the assumption it follows that $S(\Gamma)=T(\Gamma)$ and $S(\Gamma')=T(\Gamma')$. 

 In Figure~\ref{fig:triangle} the adjoint convex lattice polygon is a standard triangle of length $2$. The cornerpoints are denoted by $p_1,p_2$ and $p_3$. 
\begin{figure}[h!]
 \centering {\includegraphics[width=4cm,height=4cm]{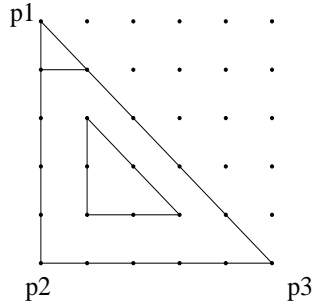}} \caption{The outer convex lattice polygon without $p_1$ is not a standard triangle, and its adjoint is a standard triangle.} 
\label{fig:triangle} 
\end{figure}

\textit{Claim 3:} If A2($\Gamma$) then $S(\Gamma)=T(\Gamma)=\{~ h \in S(\Gamma') ~~|~~ h $ tight for $\Gamma ~\}$. 
  From Lemma~\ref{lem:va_low} it follows that $\Gamma\subset$ConvexHull$(p_1,p_2,p_3)$.  At least either $p_1,p_2$ or $p_3$ is not contained by $\Gamma$, otherwise we are in case A1.  For any of these three points not contained in $\Gamma$, the direction of the opposite edge is optimal and tight.

\textit{Claim 4:} If A3($\Gamma$) then $S(\Gamma)=T(\Gamma)=S(\Gamma')$. 
  If $\Gamma'$ is not a thin triangle then it follows from Proposition~\ref{prop:va_min}.a and Lemma~\ref{lem:va_tight}.  If $\Gamma'$ is a thin triangle then it follows from Proposition~\ref{prop:va_min}.b. 

 The multiple adjoints $\Gamma^i$ for $i\in{\textbf{Z}}_{\geq 0}$ are defined in  Definition~\ref{def:adjoint_polygon}.
 We define A$n$($Z$) for $n=0,1,2,3,4$ to be as in Definition~\ref{def:casedist_surface}, but with $\Gamma$ replaced by $Z$ and $\Gamma'$ by $Z'$.
 Let \[ \alpha: {\mathcal{V}} \rightarrow {\textbf{Z}}_{\geq0},~~ \Gamma \mapsto \min_{i\geq 0}\{~ i\geq 0 {\textrm{~and~}} (\textrm{A2}(\Gamma^{i+1}) {\textrm{~or~}} \textrm{A3}(\Gamma^{i+1}))~\} \] where ${\mathcal{V}}$ is the set of all convex lattice polygons.

\textit{Claim 5:} If A4($\Gamma$) then $T(\Gamma')\neq\emptyset$ and $S(\Gamma)=S(\Gamma')$. 
  \textit{Induction claim:} $C[i]:$ If $\alpha(\Gamma)=i$ and A4($\Gamma$) then $T(\Gamma')\neq\emptyset$ and $S(\Gamma)=S(\Gamma')$, for all $\Gamma$. \textit{Induction basis $C[0]$:}  From claim 3,4) it follows that $S(\Gamma')=T(\Gamma')$.  From claim 2) is follows that $C[0]$ holds for both cases. \textit{Induction step ($C[i-1]\Rightarrow C[i]$ for $i>0$): }  We are in case A4($\Gamma'$).  From the induction hypothesis $C[i-1]$ it follows that $T(\Gamma^2)\neq\emptyset$.  From claim 1)  it follows that $S(\Gamma')=T(\Gamma')\neq\emptyset$.  From claim 2) it follows that $S(\Gamma)=S(\Gamma')$. 

 \end{proof} 
\end{theorem}

\begin{remark}\label{rem:}\textrm{\textbf{(maximal number of optimal viewangles)}} From Proposition~\ref{prop:va_min} and Figure~\ref{fig:va_opt}.16 it follows that $\#S(\Gamma)\leq 4$ for all convex lattice polygons $\Gamma\subset{\textbf{R}}^2$. Recently \cite{draisma1} proved a generalization of this result to higher dimension. They give an upperbound of $(3^d-1)/2$ for the number of optimal viewangles for the more general $d$ dimensional convex bodies in ${\textbf{R}}^d$. Moreover they show that the upperbound is only reached by the regular crosspolytopes. 
\end{remark}

\section{Families}
 
\label{sec:fam}
\begin{definition}
\label{def:fam_subset}
\textrm{\textbf{(family of subsets)}} 
  A \textit{family of subsets} $\tilde{F}$ of $\tilde{X}$ is defined as the map \[\chi: \tilde{I} \rightarrow \mathcal{P}(\tilde{X}),\quad i \mapsto \tilde{F_i},\] where  
 $\tilde{X}$ is a set, 
 $\tilde{I}$ is a set, 
 $\tilde{U}$ is a subset of $\tilde{I} \times \tilde{X}$ and 
 $\tilde{F_i}:=\{~ x \in \tilde{X} ~~|~~ (i,x) \in \tilde{U} ~\}$ for $i\in \tilde{I}$.  We defined $\tilde{F}$ to give some intuition for Definition~\ref{def:fam}. 
 
\end{definition}

\begin{definition}
\label{def:fam}
\textrm{\textbf{(family)}} 
  A \textit{family} $F$ of $X$ is defined as $(F_i)_{i\in I}$ where  
 $X$ is a projective surface over the field ${\textbf{C}}$ of complex numbers, 
 $I$ is a nonsingular curve, 
 $U$ is an irreducible, codimension 1, algebraic subset of $I \times X$ and 
 $F_i=\pi_{2*}\circ \pi_1^{-1}(\{i\})$ is an irreducible,  codimension 1, algebraic subset of $X$ for generic $i\in I$.  The maps $\pi_1: U \rightarrow I$ and $\pi_2: U \rightarrow X$ denote the first  respectively  second projection of $U$. We define $\textrm{Fam} X$ to be the set of all families on $X$.
 
\end{definition}

\begin{definition}
\label{def:}
\textrm{\textbf{(degree and geometric genus of a family)}} 
 Let $F=(F_i)_{i\in I}$ be a family as defined in Definition~\ref{def:fam}. The \textit{degree of a family}  with respect to   a given embedding $X\subset{\textbf{P}}^r$ is defined as $\deg F:=\deg F_i$ for generic $i$. The \textit{geometric genus} of a family is defined as $p_g F:=p_g F_i$ for generic $i$.
 
\end{definition}

\begin{definition}
\label{def:attr_fam}
\textrm{\textbf{(attributes of families: fibration and rational)}} 
 Let $F=(F_i)_{i\in I}$ be a family as defined in Definition~\ref{def:fam}. We call $F$ a \textit{fibration family}  if and only if  there exists a rational map  \[ f: X \dashrightarrow I \]  such that  $F_i=f^{-1}(i)^-$ for all $i\in I$.  We call $F$ a \textit{rational family}  if and only if  $p_g F=0$.
 
\end{definition}

\begin{proposition}\label{prop:fam_prop}\textbf{\textrm{(properties of families)}}   Let $F\in\textrm{Fam} X$ be a family. 
\begin{itemize} 
\item[\textbf{a)}] We have that $(F_i)_{i\in I}$ and $U$ are different representations for the same family $F$.
\item[\textbf{b)}] We have that $\textrm{supp}(F_i)=\{~ x \in X ~~|~~ (i,x) \in U ~\}$.
\item[\textbf{c)}] If $F$ is a fibration family then $\pi_2$ is birational and $f=\pi_1 \circ \pi_2^{-1}$ is a fibration map.
\item[\textbf{d)}] If $X$ is nonsingular then $F_i$ is a Cartier divisor for all $i\in I$.
\item[\textbf{e)}] If $X$ is nonsingular then $U\subset X\times I$ is a Cartier divisor.
\item[\textbf{f)}] If $X$ is nonsingular then $\deg F=\deg F_i$ for all $i\in I$ and $p_g F=\max\{~ p_g(F_i) ~~|~~ i\in I ~\}$. 
\end{itemize}
\begin{proof}  
 We have that \textbf{a)} until \textbf{e)} are straightforward.
 See \cite{har1} Corollary III.9.10 for the proof of \textbf{f)}.
 \end{proof}
\end{proposition}

\begin{proposition}\label{prop:fam_rat_prop}\textbf{\textrm{(properties of rational families)}}  Let $X$ be nonsingular.  Let $K$ be the canonical divisor class of $X$.  Let $F\in\textrm{Fam} X$ be a family. 
\begin{itemize} 
\item[] If $p_g(F)=0$ then $F K\leq-2$. 
\end{itemize}
\begin{proof}  
 Let $U\subset I\times X$ be the Cartier divisor defining $F$.
 Let $g:\tilde{U}\rightarrow U$ be the resolution of singularities of $U$ (see \cite{har1} for resolution of singularities).
 Let $\rho_1:=\pi_1\circ g: \tilde{U}\rightarrow I$ and $\rho_2:=\pi_2\circ g: \tilde{U}\rightarrow X$.
 Let $G=(G_i)_{i\in I}\in\textrm{Fam} {\tilde{U}}$ where $G_i=g^*\pi_1^{-1}(\{i\})$.

\textit{Claim 1:} We have that $G^2=0$. 
  We have that $\rho_1(G_i)\neq\rho_1(G_j)$ for all $i,j\in I$  such that  $i\neq j$.  From $\rho_1$ being a morphism it follows that $G_i\cap G_j=\emptyset$ for all $i,j\in I$  such that  $i\neq j$.

\textit{Claim 2:} If $p_g(F)=0$ then $p_a(G)=0$. 
  From $F_i=\pi_2\circ\pi_1^{-1}(\{i\})$ it follows that          $\pi_2: g(G_i) \stackrel{\cong}{\rightarrow} F_i$.  It follows that $G_i$ and $F_i$ are birational for all $i$ and thus $p_g(G)=0$.  From Sard's theorem it follows that the generic fibre $G_i$         of the regular map $\rho_1$ is nonsingular.  It follows that $p_a(G)=p_g(G)=0$. 

 Let $R=K_{\tilde{U}}-{\rho_2}^{*}K_{X}$ be the relative canonical divisor.

\textit{Claim 3:} We have that $GR\geq 0$. 
   Since we can pull back differential forms along a morphism it follows that  $0\rightarrow\rho_2^*\omega_X\rightarrow\omega_{\tilde{U}}$. From the tensor product with an invertible sheaf being exact it follows that  $0\rightarrow{\textit{O}}_X\rightarrow\omega_{\tilde{U}}\overset{}{\underset{}{\otimes}}(\rho_2^*\omega_X)^{-1}$ is exact. From the global section functor being left exact it follows that  $\omega_{\tilde{U}}\overset{}{\underset{}{\otimes}}(\rho_2^*\omega_X)^{-1}={\textit{O}}_{\tilde{U}}(R)$ is effective. From $G$ having no fixed components and being movable it follows that $G$ is nef and thus $GR\geq 0$. 

 Let (AF) denote the Adjunction Formula: $p_a(C)=\frac{1}{2}(C^2+CK)+1$ for all irreducible curves $C\subset X$ (see \cite{har1}).

\textit{Claim 4:} If $p_g(F)=0$ then $F K\leq-2$. 
  From (AF) and claim 2) it follows that $GK_{{\tilde{U}}}=2p_a(G)-G^2-2=-2$.  We have that $F K_{X}=\rho_{2*}GK_{X}=G {\rho_2}^{*}K_{X}=G K_{{\tilde{U}}}-G R\leq-2$. 

 \end{proof} 
\end{proposition}

\begin{example}\label{ex:fib_fam}\textrm{\textbf{(fibration family)}}
  Let $X={\textbf{P}}^2$.  Let $I={\textbf{P}}^1$.  Let $U=\{~ (i_0:i_1)\times(x_0:x_1:x_2) ~~|~~ x_0i_1=x_1i_0 ~\}$. 

 The corresponding family $F$ is the family of lines through a point. 

 It is a fibration family with fibration map $f: X \dashrightarrow I, ~~~~ (x_0:x_1:x_2) \longmapsto (x_0:x_1)$ 
 
\end{example}

\begin{example}\label{ex:non_fib_fam}\textrm{\textbf{(non fibration family)}}
  Let $X={\textbf{P}}^2$.  Let $I: i_0^2 + i_1^2 - i_2^2 = 0 \subset {\textbf{P}}^2$.  Let $U=\{~ (i_0:i_1:i_2)\times(x_0:x_1:x_2) ~~|~~         i_0 x_0 + i_1 x_1 - i_2 x_2 = 0 ~\}$. 

 The corresponding family $F=(F_i)_{i\in I}$ is the family of tangents to a circle in a plane. 

 The family $F$ is not a fibration family. 

 The intersection of two lines is varying with the pair of lines. In other words, generic points in $X$ are reached by $2$  family members $F_i$. 
 
\end{example}

\begin{definition}
\label{def:}
\textrm{\textbf{(operations on families)}} 
  Let $F\in\textrm{Fam} X$ as in Definition~\ref{def:fam}. Let $f: X \rightarrow Y$ be a birational morphism between projective surfaces. The \textit{pushforward of families} is defined as \[ \ensuremath{f_\circledast:\textrm{Fam} X\rightarrow \textrm{Fam} Y,\quad U\mapsto \hat{f}(U)}. \] The \textit{pullback of families} is defined as \[ \ensuremath{f^\circledast:\textrm{Fam} Y\rightarrow \textrm{Fam} X,\quad V\mapsto \overline{\hat{f}^{-1}(V-B)}}. \] where  
 $\ensuremath{\hat{f}:I\times X\rightarrow I\times Y,\quad (~i,~x~)\mapsto (~i,~f(x)~)}$ and 
 $B\subset I\times Y$ is the locus where $\hat{f}^{-1}$ is not defined.  If $X$ is nonsingular then the intersection products are defined as  
 $\cdot: \textrm{Div} X \times \textrm{Fam} X \rightarrow {\textbf{Z}},~~ (D,F) \mapsto DF_i \textrm{~for any~} i\in I$ and 
 $\cdot: \textrm{Fam} X \times \textrm{Fam} X \rightarrow {\textbf{Z}},~~ (F,F') \mapsto F_jF_i' \textrm{~for any~} i\in I {\textrm{~and~}} j\in I'$.  The following proposition shows that the intersection products are well defined. 

\end{definition}

\begin{proposition}\label{prop:}\textbf{\textrm{(properties of  operations on families)}}  Let $h:X\rightarrow Y$ be a birational morphism between surfaces. 
\begin{itemize} 
\item[\textbf{a)}] The maps $h_\circledast$ and $h^\circledast$ are well defined.
\item[\textbf{b)}] We have that $h^\circledast \circ h_\circledast=\textrm{id}_{\textrm{Fam} X}$ and $h_\circledast \circ h^\circledast=\textrm{id}_{\textrm{Fam} Y}$.
\item[\textbf{c)}] If $X$ and $Y$ are nonsingular then   
 $\ensuremath{f_\circledast:\textrm{Fam} X\rightarrow \textrm{Fam} Y,\quad (F_i)_{i\in I}\mapsto }  {(f_*F_i)_{i\in I}},~~$ and 
 $\ensuremath{f^\circledast:\textrm{Fam} Y\rightarrow \textrm{Fam} X,\quad (F_i)_{i\in I}\mapsto } {\overline{(~~(f^* F_i)_{i\in I}-(\cap_{i\in I}F_i)~~)}}$  where $f^*$ and $f_*$ are defined by the pullback and pushforward of divisors.
\item[\textbf{d)}] If $X$ is nonsingular then $DF_i=DF_j$ for all $D\in\textrm{Div} X$ and $i,j\in I$ and thus the intersection products are well defined. 
\end{itemize}
\begin{proof} 
 We have that \textbf{a)}, \textbf{b)} and \textbf{c)} are a straightforward consequence of the definitions.
 See \cite{har1} for the proof of \textbf{d)} (family members $F_i$ are algebraic equivalent and algebraic equivalence implies numerical equivalence).
 \end{proof}
\end{proposition}

\section{Minimal degree families on toric surfaces}
 
\label{sec:fam_on_tor}
\begin{remark}\label{rem:tor_var}\textrm{\textbf{(toric varieties)}} For the definition of toric varieties we follow  \cite{ewa1}, \cite{cox1} and \cite{ful1}. If $\Gamma$ is a lattice polygon with lattice points $\{(a_0,b_0),\dots,(a_r,b_r)\}$, then the toric surface defined by $\Gamma$  is the projective closure of the image of the map \[ p:{\textbf{C}}^{*2}\rightarrow{\textbf{P}}^r, (s,t)\mapsto (s^{a_0}t^{b_0}:\dots:s^{a_r}t^{b_r}) \] (see \cite{cox1} section 12). 
\end{remark}

\begin{definition}
\label{def:}
\textrm{\textbf{(attributes of families: toric family)}} 
  Let $F$ in $\textrm{Fam} X$ be a family as defined in Definition~\ref{def:fam}.  We call $F$ a \textit{toric family}  if and only if  $F$ is a fibration family and after resolution of basepoints the fibration map is a toric morphism. Note that $X$ and $I$ have to be toric and in particular $I={\textbf{P}}^1$ (see \cite{ewa1} for the definition of toric morphism).  The fibration map induces a toric morphism between the dense tori in $X$ and $I$ (see Example~\ref{ex:tor_va_rel} below).  
 
\end{definition}

\begin{definition}
\label{def:torfam}
\textrm{\textbf{(minimal toric degree and optimal toric family)}}  Let $X$ be a complex embedded toric surface. 
 The \textit{minimal toric degree} $v(X)$ of $X$  is the smallest possible degree of a toric family on $X$ (see Definition~\ref{def:attr_fam}). The \textit{set of optimal toric families} on $X$ is defined as  \[ S(X)= \{~ F\in \textrm{Fam} X ~~|~~ F \textrm{ is a toric family }  {\textrm{~and~}}   \deg F = v(X) ~\}. \]
 
\end{definition}

\begin{definition}
\label{def:toric_problem}
\textrm{\textbf{(toric family problem on toric surfaces)}} 
  Given a complex embedded toric surface $X$ find the minimal toric degree $v(X)$ and the set of optimal toric families $S(X)$. 
 
\end{definition}

\begin{definition}
\label{def:tor_va_rel}
\textrm{\textbf{(viewangles and toric families relation)}} 
  Let $\Gamma$ be a lattice polygon with lattice $\Lambda$. Let $X$ be the toric surface defined by $\Gamma$ (see Remark~\ref{rem:tor_var}). Let $V$ be the set of primitive viewangles in $\Lambda^*-\{0\}$. Let $T$ be the set of toric families on $X$. The \textit{viewangles and toric families relation} is  a function: \[ \theta_\Gamma: V \rightarrow T \] where any primitive viewangle $h\in V$ is send to a toric family in $\theta_\Gamma(h) \in T$ in the following way: 
 Let $\Sigma$ with lattice $\Lambda^*$  be the normal fan of $\Gamma$ (see \cite{cox1} section 12).
 Let $\Sigma'$ be the fan of ${\textbf{P}}^1$   (the unique projective toric curve) with lattice points in $\Lambda^*/h$ .
 Let $\tau$ and $\tau'$ be the cones in $\Sigma$  respectively  $\Sigma'$  corresponding to the dense torus embeddings (thus the cones are points).
 The canonical linear map  $\Lambda^*\rightarrow\Lambda^*/h$ induces map of fans $\alpha:\tau\rightarrow\tau'$  (see \cite{ewa1} section V.4 for map of fans).
 Let $\beta:X_\tau\rightarrow X_{\tau'}$ be the toric morphism corresponding to the map of fans $\alpha$ (see \cite{ewa1} section VI.6). 
 Let $f:X_{\Sigma}\dashrightarrow X_{\Sigma'}$ be the rational map corresponding to the closure of $\beta$.
 The toric family $\theta_\Gamma(h)$ is defined by the fibres of $f$. 
 
\end{definition}

\begin{example}\label{ex:tor_va_rel}\textrm{\textbf{(viewangles and toric families relation)}} 
  Let $\theta_\Gamma: V \rightarrow T$ be  the viewangles and toric families relation. We use the same notation as in Definition~\ref{def:tor_va_rel}.

 We assume that $\Gamma$ with lattice $\Lambda$ is the standard triangle in Figure~\ref{fig:exm_tor_relQ}.a). The vertical lines represent the viewangle $h=(m,n)=(0,-1)$ in $V$. 
\begin{figure}[h!]
 \centering \begin{tabular}
{|c|c|c|c|c|} \hline $\Gamma$ and $h$ & $\tau\stackrel{\alpha}{\rightarrow}\tau'$  \\
\hline & \\
 {\includegraphics[width=2cm,height=2cm]{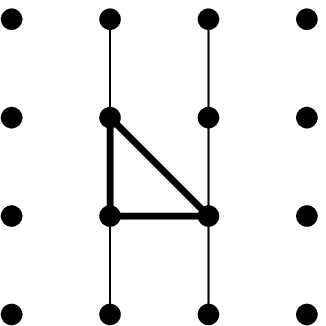}}&{\includegraphics[width=2cm,height=2cm]{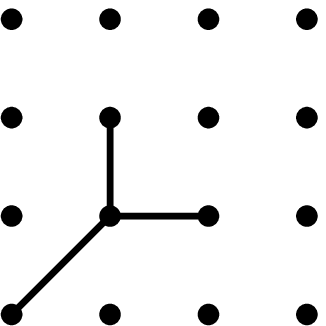}} \\
 & $\downarrow\alpha$ \\
 &{\includegraphics[width=2cm,height=0.1cm]{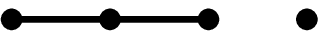}} \\
\hline a) & b) \\
\hline 
\end{tabular}
 \caption{Example of toric families and viewangles relation.} 
\label{fig:exm_tor_relQ} 
\end{figure}

 The triangle polytope in Figure~\ref{fig:exm_tor_relQ}.a) corresponds to the closure of  the image of \[ p:{\textbf{C}}^{*2}\rightarrow{\textbf{P}}^2, (s,t)\mapsto (s:t:1)  \] which is ${\textbf{P}}^2$. 

 In Figure~\ref{fig:exm_tor_relQ}.b) is the normal fan $\Sigma$ of the triangle polygon with lattice $\Lambda^*$. Downstairs is the fan of ${\textbf{P}}^1$ which is the unique projective toric curve, with lattice $\Lambda^*/h$.

 The canonical linear map $\Lambda^*\rightarrow\Lambda^*/h$ is defined by the matrix $[n~~ m]=[-1~~ 0]$, which is the vertical projection. 

 It induces a map of fans $\beta:\tau\rightarrow\tau'$ on the dense torus embeddings (see Figure~\ref{fig:exm_tor_relQ}.b)).

 The map $\beta$ defines a semigroup homomorphism: \[ \beta^*: [u,u^{-1}] \rightarrow [s,t,s^{-1},t^{-1}],\quad u \mapsto s^{-n} t^{m}.   \]

 We have that $\beta^*$ defines the following rational map between  the toric varieties: \[ f': {\textbf{C}}^{*2} \rightarrow {\textbf{C}}^{*},\quad (s,t) \mapsto (s^{-n} t^{m}) \] The closure of $f'$ defines the map \[ f: {\textbf{P}}^2\dashrightarrow {\textbf{P}}^1,\quad (x_0:x_1:x_2) \mapsto  (x_0:x_2)=(x_0^{-n} x_1^{m} x_2^{n-m}:1) \] which is not defined at $(0:1:0)$. 

 The corresponding toric family $\theta_\Gamma(h)$  is the family of lines through the  point $(0:1:0)$. 

 This family has degree $1$ and $h$ is an optimal viewangle  of width $1$ (see Figure~\ref{fig:va_opt}.17). This is no coincidence  as we shall see in Proposition~\ref{prop:tor_va_rel}.
 
\end{example}

\begin{proposition}\label{prop:tor_va_rel}\textbf{\textrm{(viewangles and toric families relation)}}  Let $\theta_\Gamma: V \rightarrow T$ be  the viewangles and toric families relation. 
\begin{itemize} 
\item[\textbf{a)}] We have that $\theta_\Gamma$ is a bijection and a viewangle of width $n$ is send to a toric family of degree $n$. 
\end{itemize}
\begin{proof} 
 We use the same notation as in Definition~\ref{def:tor_va_rel}.
  Let $\{(a_0,b_0),\dots,(a_r,b_r)\}$ be the set of lattice  points of $\Gamma$.
 Let $p: {\textbf{C}}^{*2} \rightarrow {\textbf{P}}^r,\quad (s,t) \mapsto (s^{a_0}t^{b_0}:\dots:s^{a_r}t^{b_r})$ (see Remark~\ref{rem:tor_var}).
 Let $h=(m,n)$ in $V$ be a primitive viewangle ( id est  $\gcd(m,n)=1$).
 Let $F=\theta_\Gamma(h)$.
 Let \[ f: X \dashrightarrow {\textbf{P}}^1,\quad (x_0:\dots:x_r) \mapsto ({x_0}^{e_0}\dots{x_r}^{e_r}:1) \]  such that  $\sum e_i=0$, $\sum a_ie_i=-n$ and  $\sum b_ie_i=m$ for $i\in\{0,\ldots,r\}$.

\textit{Claim 1:} We have that $f$ is the fibration map  of toric family $F$.  
  In Example~\ref{ex:tor_va_rel} the map $f$ is obtained for a special case. That this construction holds in general is left to the reader. 

 Let $q=(q_0:q_1)\in{\textbf{P}}^1$.

\textit{Claim 2:} The fibres $f(q)^{-1}$ are  $F_q:=\overline{\{~x \in X ~|~ {x_0}^{e_0}\dots{x_n}^{e_n} = \frac{q_0}{q_1}~\}}$. 
  This claim is a direct consequence of the definitions.

\textit{Claim 3:} We have that $p^{-1}(F_q): s^{-n}t^m - \frac{q_0}{q_1}=0$ and  this curve is irreducible  if and only if  $\gcd(m,n)=1$. 
  If $\alpha,\beta$ are coprime and $z\in{\textbf{Z}}_{>1}$ then $s^{z\alpha}t^{z\beta}-1=(s^{\alpha}t^{\beta}-1)(\sum_{i=0}^{i=z}s^{i\alpha}t^{i\beta})$. 

 Let $k,l \in {\textbf{C}}^*$ be  such that  $\frac{q_0}{q_1}=\frac{l^m}{k^n}$.
 Let $h_q: {\textbf{C}}^* \to {\textbf{C}}^{*2},\quad u \mapsto (k \cdot u^m, l \cdot u^n)$.

\textit{Claim 4:} The map $h_q$ is a birational parametrization of $p^{-1}(F_q)$. 
  This claim is a direct consequence of the definitions. 

 Let $g_{mn}(q): {\textbf{C}}^* \to F_q,\quad u \mapsto (k^{a_0} l^{b_0} \cdot u^{a_0 m + b_0 n}:\dots:k^{a_r} l^{b_r} \cdot u^{a_r m + b_r n})$.

\textit{Claim 5:} The map $g_{mn}(q)$ is a birational parametrization of $F_q$ for all generic $q\in{\textbf{P}}^1$. 
  We have that $g_{mn}(q)=p\circ h_q$ for all $q\in{\textbf{P}}^1$.  We have that $f\circ g_{mn}(u)=k^{\sum a_i e_i}l^{\sum b_i e_i}u^{\sum a_i e_i m + \sum b_i e_i n}$.  It follows that $\sum a_i e_i=-n$ and $\sum b_ie_i=m$.

\textit{Claim 6:} Changing $k,l$ in $g_{mn}(q)$  such that  $\frac{q_0}{q_1}=\frac{l^m}{k^n}$  gives rise to a reparameterization of $F_q$. 
  Direct consequence of the definition of $h_q$ and  that $g_{mn}(q)=p\circ h_q$ for all $q\in{\textbf{P}}^1$.

\textit{Claim 7:} We have that  $\deg F=\displaystyle\max_i(a_i m+b_i n) - \displaystyle\min_i(a_i m+b_i n)$. 
  From claim 5) it follows that $\deg(F)$ equals the cardinality of $g_{mn}(q)\cap H$ for any $q$ and generic hyperplane section $H$.

\textit{Claim 8:} We have that \textbf{a)}. 
  The linear system of equations $\sum a_i e_i=-n, \sum b_ie_i=m$ and $\sum e_i=0$ has solutions in $\vec{e}$.  From claim 6) it follows that  $F$ corresponding to $g_{m'n'}$ depends uniquely on $a_i,b_i,m'$ and $n'$.  It follows that $\theta_\Gamma(h)$ defines uniquely a family $F$.   From claim 7) it follows that a viewangle of width $n$ is send  to a toric family of degree $n$. 

 \end{proof} 
\end{proposition}

\begin{example}\label{ex:toricfamily}\textrm{\textbf{(toric family problem)}} 
  Let $X$ be a complex embedded toric surface.  Let $p: {\textbf{C}}^{*2} \to X \subset {\textbf{P}}^6,\quad (s,t) \mapsto (s^0 t^1:s^0 t^2:s^1 t^0:s^1 t^1:s^1 t^2:s^2 t^0:s^2 t^1)$ a birational monomial parameterization.

 For $m=1$ and $n=-1$ we find \[ f: X \dashrightarrow {\textbf{P}}^1,\quad (x_0:\dots:x_6) \mapsto (x_1^2 x_2:x_0^3) \] and $\deg(F_q)=4$ for all $q\in{\textbf{P}}^1$. 

 For $m=1$ and $n=0$ we find \[ f: X \dashrightarrow {\textbf{P}}^1,\quad (x_0:\dots:x_6) \mapsto (x_1:x_0) \] and $\deg(F_q)=2$ for all $q\in{\textbf{P}}^1$. We have that Example~\ref{ex:toricfamily} and Example~\ref{ex:va_problem} reflect an equivalent problem instance.  
 
\end{example}

\section{Adjoint chain}
 
\label{sec:adjoint_chain}
\begin{remark}\label{rem:}\textrm{\textbf{(references)}} We claim no new results in this section. For the notion of nef, movable, canonical class and exceptional curve we refer to \cite{har1} and \cite{mat1}. The adjoint chain is a reformulation and adapted version of $(D+K)$-minimalization as described in \cite{man1} and can also be found  in \cite{Schicho:97}. 
\end{remark}

\begin{definition}
\label{def:mprs}
\textrm{\textbf{(minimally polarized rational surface (mprs))}} 
  A \textit{minimally polarized rational surface (mprs)} is  defined as a pair $(X,D)$ where  
 $X$ is a nonsingular rational surface over ${\textbf{C}}$,  
 $D$ is a nef and movable divisor on $X$ and 
 there doesn't exists a $-1$-curve $C$  such that  $DC=0$.  
 
\end{definition}

\begin{definition}
\label{def:}
\textrm{\textbf{(minimal mprs)}} 
  Let $(X,D)$ be a mprs.  Let $K$ denote the canonical divisor class on $X$. We call $(X,D)$ a \textit{minimal mprs}   if and only if  $\dim|D+K|\leq 0$ or $D^2=0$. 
 
\end{definition}

\begin{definition}
\label{def:adjoint_relation}
\textrm{\textbf{(adjoint relation)}} 
  Let $(X,D)$ be a mprs which is not minimal. An \textit{adjoint relation} is a relation  $(X,D) \stackrel{\mu}{\rightarrow} (X',D')$ where  
 $(X,D)$ is a mprs which is not minimal, 
 $(X',D')$ is a mprs, 
 $\ensuremath{X\stackrel{\mu}{\rightarrow}X'}$ is a birational morphism  which blows down all $-1$-curves $C$  such that  $(D+K)C=0$ and 
 $D'=\mu_*(D+K)$. 
 
\end{definition}

\begin{definition}
\label{def:adjoint_chain}
\textrm{\textbf{(adjoint chain)}} 
  An \textit{adjoint chain} of $(X,D)$ is a chain of adjoint relations  until a minimal mprs is obtained: \[ \ensuremath{(X,D)=(X_0,D_0)\stackrel{\mu_0}{\rightarrow}(X_1,D_1)} \ensuremath{\stackrel{\mu_1}{\rightarrow}} \ldots. \]  
 
\end{definition}

\begin{proposition}\label{prop:adjoint_chain}\textbf{\textrm{(properties of adjoint chain)}}
\begin{itemize} 
\item[\textbf{a)}] The adjoint chains of a mprs are finite and have the same length. 
\item[\textbf{b)}] If $\ensuremath{(X,D)\stackrel{\mu}{\rightarrow}(X',D')}$ is an adjoint relation then $\mu^* D' = D+K$. 
\end{itemize}
\begin{proof}  
 The proofs can be found in \cite{Schicho:97}. 
 \end{proof} 
\end{proposition}

\section{Minimal degree families on polarized rational surfaces}
 
\label{sec:fam_on_mprs}
\begin{definition}
\label{def:fam_mprs}
\textrm{\textbf{(optimal and tight families and minimal degree)}} 
  Let $(X,D)$ be a mprs.  Let $K$ be the canonical divisor class on $X$.  Let $F\in\textrm{Fam}(X)$. The \textit{degree} of $F$  with respect to  $(X,D)$ is given by $DF$. We call $F$ a \textit{tight family}  if and only if  $FK=-2$. The \textit{minimal rational degree}  with respect to  $(X,D)$ is defined as \[ v(X,D)=\min\{~ DF ~~|~~ F\in\textrm{Fam} X {\textrm{~and~}} p_g(F)=0 ~\}. \] The minimum exists since $D$ is nef by definition.  We call $F$ an \textit{optimal family}  if and only if  $F$ is a rational family and $DF=v(X,D)$. The \textit{set of all optimal families} on $(X,D)$ is denoted by $S(X,D)$. 
 
\end{definition}

\begin{example}\label{ex:}\textrm{\textbf{(optimal families of the projective plane)}} 
  Let $F$ be the family of lines through a point (see Example~\ref{ex:fib_fam}).  Let $L$ be the divisor class of lines on ${\textbf{P}}^2$. 

 We have that $({\textbf{P}}^2,L)$ is a mprs. 

 We have that $F\in S({\textbf{P}}^2,L)$ and $v({\textbf{P}}^2,L)=FL=1$. 
 
\end{example}

\begin{definition}
\label{def:fam_mprs_problem}
\textrm{\textbf{(rational family problem on mprs)}} 
  Given a mprs $(X,D)$ find the minimal degree $v(X,D)$ and all optimal families $S(X,D)$. 
 
\end{definition}

\begin{lemma}\label{lem:fam_low}\textrm{\textbf{(lowerbound)}}  Let $\ensuremath{(X,D)\stackrel{\mu}{\rightarrow}(X',D')}$ be an adjoint relation.  Let $F\in\textrm{Fam} X$ be a rational family. 
\begin{itemize} 
\item[] We have that $F D\ge \mu_\circledast FD' + 2$, and equality holds   if and only if   $F$ is tight. 
\end{itemize}
\begin{proof} 
 From Proposition~\ref{prop:adjoint_chain}.b and Proposition~\ref{prop:fam_rat_prop} it follows that $\mu_\circledast F D'=F \mu^*D'=F D+ F K \leq F D - 2$.
 \end{proof} 
\end{lemma}

\begin{lemma}\label{lem:fam_tight}\textrm{\textbf{(tight)}}  Let $\mu: X\rightarrow X'$ a birational morphism between nonsingular complex projective surfaces.  Let $F'\in\textrm{Fam} X'$ be tight. 
\begin{itemize} 
\item[\textbf{a)}] We have that $\mu^\circledast F'\in\textrm{Fam} X$ is tight. 
\item[\textbf{b)}] We have that $\mu^\circledast F'=(\mu^* F'_i)_{i\in I}$. 
\end{itemize}
\begin{proof}  

\textit{Claim:} We assume  without loss of generality  that $\mu=\pi$ where $\pi$ blows down one exceptional curve $E$. 

\textit{Claim:} We have that \textbf{a)} and \textbf{b)}. 
  We have that $\mu^\circledast F' K= F' \mu_* K = F' K' = -2$.  We have that $\mu^\circledast F' K = (F + mE) K = F K - m=-2$ where $m\ge 0$.  From Proposition~\ref{prop:fam_rat_prop} it follows that $m=0$ and $F K=-2$. 
 
 \end{proof}
\end{lemma}
  
\begin{proposition}\label{prop:fam_min}\textbf{\textrm{(classification optimal fibration families on minimal mprs)}} 
\begin{itemize} 
\item[] All the optimal \textit{fibration} families on minimal mprs $(X,D)$ are classified in the following table: \\[2mm]
 {\tiny \begin{tabular}
{|c|c|c|c|c|c|c|} \hline $D$            & $D^2$ & $X\cong{\textbf{P}}^2$ & optimal families $F$      & $DF$ & tight & type \\
\hline $D=nP$         & $0$   & no             & $F=P$                     & $0$  & yes   & ruled \\[2mm]
 $2D+K=nP$      & $n+2$ & no             & $F=P$                     & $1$  & yes   & linear fibration \\[2mm]
 $2D+K=0$       & $2$   & no             & 1 or 2 families of lines  & $1$  & yes   & linear fibration \\[2mm]
 $D=L$          & $1$   & yes            & $F\subset L$              & $1$  & no    & plane \\[2mm]
 $D=2L$         & $2$   & yes            & $F\subset L$              & $2$  & no    & plane \\[2mm]
 $D+K=0$        & $9$   & yes            & $F\subset L$              & $3$  & no    & plane \\[2mm]
 $D+K=0$&$1,2,\ldots,8$ & no             & see \cite{Schicho:00a}    & $2$  & yes   & conic fibration \\[2mm]
 $D+cK=0$       & $0$   & no             & infinitely many           & $2c$ & yes   & \\
\hline 
\end{tabular}
 } \\[2mm]
 where  
 $p_g(L)=0$,~$\ensuremath{\textrm{dim}}|L|=2$ and $L^2=1$ ($L$ stands for lines); and 
 $p_g(P)=0$,~$\ensuremath{\textrm{dim}}|P|=1$ and $P^2=0$ and 
 $c\in{\textbf{Z}}_{>0}$.  In particular we see that there is always an optimal family of fibration type. 
\end{itemize}
\begin{proof}  
 The first 3 columns are known from \cite{man1}.
 The third row denotes families of lines of a quadric surface in ${\textbf{P}}^3$.
 The rows 4 to 7 are known from \cite{Schicho:00a} (page 81 until 85). The cases $D^2=1,2$ in row 7 are not covered in \cite{Schicho:00a}, but are straightforward generalizations.
 The last row is the Halphen pencil and can be found in \cite{hal1} and Exercise V.4.15.e in \cite{har1}. This pair can never arise as a last link in an adjoint chain where the mprs $(X_0,D_0)$ satisfies $D_0^2>0$.
 Let (AF) denote the Adjunction Formula: $p_a(C)=\frac{1}{2}(C^2+CK)+1$ for all irreducible curves $C\subset X$ (see \cite{har1}).
 Let $F=(F_i)_{i\in I}$ in $\textrm{Fam} X$ be any family  such that  $FP=0$.

\textit{Claim 1:} We have that $F=P$. 
  From $FP=0$ and $F,P$ being movable it follows that there exist curves $C\in|P|$ and $F_j \in F$ through some generic point $x\in X$.  From $CF_j=0$ and $x\in C\cap F_j$ it follows that $C=F_j$ and thus $F=P$.

\textit{Claim 2:} If $D=nP$ then $P$ is the unique optimal tight fibration family. 
  From (AF) it follows that $p_a(P)=\frac{1}{2}(0+PK)+1=0$, and thus $PK=-2$.  From $D$ being nef and $DP=0$ it follows that $P$ is an optimal family.  The fibration map is given by $\varphi_{|D|}$.  From claim 1) it follows that $P$ is the unique optimal family.

\textit{Claim 3:} If $2D+K=nP$ then $P$ is the unique optimal tight fibration family. 
  We have that $F(2D+K)\geq 0$ for all $F\in S(X,D)$.  From Proposition~\ref{prop:fam_rat_prop} it follows that $2FD=F(2D+K)-FK\geq 0 +2$  and thus $FD\geq 1$.  If $F=P$ then $FD=1$.  If $FD=1$ then $2FD=F(2D+K)-FK=2$ and thus $FP=0$.  From claim 1) it follows that $P$ is the unique optimal family. 

 \end{proof} 
\end{proposition}

\begin{definition}
\label{def:casedist_surface}
\textrm{\textbf{(case distinction)}}  Let $\ensuremath{(X,D)\stackrel{\mu}{\rightarrow}(X',D')}$ be an adjoint relation. 
  We distinguish the following cases where $(X,D)$ is not minimal except at B0: \\[2mm]
 \begin{tabular}
{|l|l|l|} \hline    & $(X,D)$          & $(X',D')$ \\
\hline B0 & minimal mprs     & $-$ \\
\hline B1 & $X\cong{\textbf{P}}^2$  & $X'\cong{\textbf{P}}^2$ \\
\hline B2 & $X\ncong{\textbf{P}}^2$ & $X'\cong{\textbf{P}}^2$ \\
\hline B3 & $X\ncong{\textbf{P}}^2$ & minimal mprs and $X'\ncong{\textbf{P}}^2$ \\
\hline B4 & $X\ncong{\textbf{P}}^2$ & not minimal mprs and $X'\ncong{\textbf{P}}^2$ \\
\hline 
\end{tabular}
 \\

\end{definition}

\begin{theorem}\label{thm:optfam}\textbf{\textrm{(optimal families and minimal degree)}}  Let $\ensuremath{(X,D)\stackrel{\mu}{\rightarrow}(X',D')}$ be an adjoint relation.  Let B0 until B4 denote the cases as in Definition~\ref{def:casedist_surface}.  Let $L$ be the divisor class of lines on $X$, if $X\cong{\textbf{P}}^2$.  Let $L'_p$ be the family of lines through the point $p$ for any $p\in X'$, if $X'\cong{\textbf{P}}^2$.  Let $B$ be the set of indeterminacy points of $\mu^{-1}$. 
\begin{itemize}
\item[\textbf{a)}] If B0 then  
 $S(\Gamma)$ and $v(\Gamma)$ are given by Proposition~\ref{prop:fam_min}. 
\item[\textbf{b)}] If B1 then  
 $S(X,D)=\{~ F ~~|~~ F\subset L ~\}$ and 
 $v(X,D) = v(X',D') + 3$. 
\item[\textbf{c)}] If B2 then  
 $S(X,D)=\{~ \mu^\circledast L'_p ~~|~~ p \in B ~\}$ and 
 $v(X,D) = v(X',D') + 2$. 
\item[\textbf{d)}] If B3 or B4 then  
 $S(X,D)=\{~ \mu^\circledast F' ~~|~~ F' \in S(X',D') ~\}$ and 
 $v(X,D)=v(X',D') + 2$. 
\end{itemize}
\begin{proof}  
 We have that \textbf{a)} is a direct consequence of Proposition~\ref{prop:fam_min}.
 We have that \textbf{b)} follows from claim 1), \textbf{c)} follows from claim 5) and \textbf{d)} follows from claim 8) and claim 9), where the claims are given below.
 Let $L$ and $L'$ be the class of lines on  respectively  $X$ and $X'$, if $X\cong{\textbf{P}}^2$ or $X'\cong{\textbf{P}}^2$.

\textit{Claim 1:} If B1 then $F\subset L$ and $v(X,D)=v(X',D')+3$. 
  If $X\cong{\textbf{P}}^2$ then $F\subset L$ for all $F\in S(X,D)$.  From Lemma~\ref{lem:fam_tight}.b and $K_{{\textbf{P}}^2}=-3L$ it follows that $L'D'=\mu^* L' \mu^* D' =L(D+K)=LD-3$. 

 Let $\ensuremath{({\tilde{X}},{\tilde{D}})\stackrel{g}{\rightarrow}(X',D')}$ be a relation  such that   
  $g:{\tilde{X}} \rightarrow X'$ is the blowup of a point $p\in B$ and ${\tilde{D}}=g^*D'$. 
 Let $\ensuremath{(X,D)\stackrel{f}{\rightarrow}({\tilde{X}},{\tilde{D}})}$ be a relation  such that   
 $\mu=g \circ f$ and ${\tilde{D}}=f_*(D+K)$. 

\textit{Claim 2:} The relation  $\ensuremath{(X,D)\stackrel{f}{\rightarrow}({\tilde{X}},{\tilde{D}})}\ensuremath{\stackrel{g}{\rightarrow}(X',D')}$  where $\mu=g \circ f$ exists. 
  It follows from \cite{har1}, Proposition V.5.3 (factorization of birational morphisms). 

 Let $G'_p=L'_p$ and  ${\tilde{G}}_p=g^\circledast G'_p= (g^*G'_{pi})_{i\in I}-{\tilde{E}}$ and  $G_p=f^\circledast{\tilde{G}}_p=\mu^\circledast G'_p$ for $p\in B$ which is blown up by $g$.
 Let $R(X,D)=\{~ G_p ~~|~~ p \in B ~\}$.
 Let $T(X,D) = \{~ F ~~|~~ F \in S(X,D) {\textrm{~and~}} F \textrm{ is tight } ~\}$.

\textit{Claim 3:} If B2 then $R(X,D)\subseteq T(X,D)$. 
  We have that ${\tilde{G}}_p{\tilde{K}}=((g^*G'_{pi})_{i\in I}-{\tilde{E}}){\tilde{K}}=G'_p g_*{\tilde{K}} + 1=-2$.  From Lemma~\ref{lem:fam_tight}.b it follows that $(f^*{\tilde{G}}_{pi})_{i\in I}=f^\circledast{\tilde{G}}_p=G_p$.  It follows that $G_pK=f^\circledast{\tilde{G}}_p K={\tilde{G}}_p f_*K=-2$.  We have that ${\tilde{G}}_p{\tilde{D}}=(g^*G'_{pi}-{\tilde{E}}){\tilde{D}}=g^*G'_{pi} g^*D'=G'_pD'$ for all $i\in I$.  From Lemma~\ref{lem:fam_low} it follows that $G_pD=f_\circledast G_p{\tilde{D}} + 2={\tilde{G}}_p{\tilde{D}} + 2=G'_pD'+2$.  From claim 1) it follows that $G'_pD'$ is minimal and thus $G_p\in T(X,D)$.

\textit{Claim 4:} If B2 then $R(X,D)\supseteq S(X,D)$. 
  If $F\in S(X,D)$ then $FD\geq\mu_\circledast FD'+2$ and thus $\mu_\circledast F \subset L'$.  It follows that $FD=\mu^\circledast\mu_\circledast FD=L'D'+2$ where $(\mu^*\mu_* F_i)_{i\in I}K=-3$.  From Lemma~\ref{lem:fam_low} it follows that $(\mu^*\mu_* F_i)_{i\in I} D > L'D'+2$.  It follows that $\mu^\circledast\mu_\circledast F\neq (\mu^*\mu_* F_i)_{i\in I}$ and they differ by a fixed component, which can only come from $p\in B$.

\textit{Claim 5:} If B2 then $S(X,D)=T(X,D)=\{~ \mu^\circledast L'_p ~~|~~ p \in B ~\}$. 
  It follows from claim 3) and claim 4).

\textit{Claim 6:} If $T(X',D')\neq\emptyset$ then $S(X,D)=T(X,D)$. 
  From Lemma~\ref{lem:fam_low} it follows that if $F'\in T(X',D')$ then $\mu^\circledast F'D = v(X',D')+2$.  From Lemma~\ref{lem:fam_low} it follows that if $F\in S(X,D)$ then $FD\geq v(X',D')+2$ and equality holds  if and only if  $F\in T(X,D)$. 

 Let $\mu^\circledast S(X,D) = \{~ \mu^\circledast F ~~|~~ F \in S(X,D) ~\}$.

\textit{Claim 7:} If $S(X',D')=T(X',D')$ then $S(X,D)=\mu^\circledast S(X',D')$. 
  From Lemma~\ref{lem:fam_tight} it follows that if $F'\in T(X',D')$ then $\mu^\circledast F'D = v(X',D')+2$.  It follows that $T(X,D)\supseteq\mu^\circledast T(X',D')$.  If $F\in T(X,D)$ then $FD=\mu_\circledast FD+2$ and thus $\mu_\circledast F\in S(X',D')$.  It follows that $T(X,D)\subseteq\mu^\circledast S(X',D')$.  From claim 6) and the assumption it follows that $S(X,D)=T(X,D)$ and $S(X',D')=T(X',D')$.

\textit{Claim 8:} If B3 then $S(X,D)=T(X,D)=\mu^\circledast S(X',D')$. 
  It follows from Proposition~\ref{prop:fam_min}, claim 5) and claim 6). 

 We will use the adjoint chain (see Definition~\ref{def:adjoint_chain}) and define $(X,D)$ to be $(X_0,D_0)$ and $(X',D')$ to be $(X_1,D_1)$.
 We define B$n$($X_i,D_i$) for $n=0,1,2,3,4$ to be as in Definition~\ref{def:casedist_surface}, but with $(X,D)$ replaced by $(X_i,D_i)$ and $(X',D')$ replaced by $(X_{i+1},D_{i+1})$.
 Let \[ \alpha: {\mathcal{V}} \rightarrow {\textbf{Z}}_{\geq0},~~ (X,D) \mapsto \min_{i\geq 0} \{~ i ~~|~~  i\geq 0 {\textrm{~and~}} (\textrm{B2}(X_{i+1},D_{i+1})  {\textrm{~or~}} \textrm{B3}(X_{i+1},D_{i+1})) ~\} \] where ${\mathcal{V}}$ is the set of all mprs's . It follows from Proposition~\ref{prop:adjoint_chain} that the length of an adjoint chain of $(X,D)$  is unique, and thus $\alpha$ is well defined.

\textit{Claim 9:} If B4$(X,D)$ then $T(X',D')\neq\emptyset$ and $S(X,D)=\mu^\circledast S(X',D')$. 
  \textit{Induction claim:} $C[i]:$ If $\alpha(X,D)=i$ and B4$(X,D)$ then $T(X',D')\neq\emptyset$ and $S(X,D)=\mu^\circledast S(X',D')$, for all $(X,D)$. \textit{Induction basis $C[0]$:}  From claim 5,8) it follows that $S(X',D')=T(X',D')$.  From claim 7) it follows that $C[0]$ holds for both cases. \textit{Induction step ($C[i-1]\Rightarrow C[i]$ for $i>0$): }  We are in case B4$(X_1,D_1)$.  From the induction hypothesis $C[i-1]$ it follows that $T(X_2,D_2)\neq\emptyset$.  From claim 6)  it follows that $S(X',D')=T(X',D')\neq\emptyset$.  From claim 2) it follows that $S(X,D)=\mu^\circledast S(X',D')$. 

 \end{proof} 
\end{theorem}

\begin{remark}\label{rem:analogy}\textrm{\textbf{(analogy with finding optimal viewangles on vineyards)}} The analogy between this section and \textsection\ref{sec:va_on_vy}, is stated in the following table: 
\begin{center}
 {\tiny \begin{tabular}
{|l|l|l|} \hline \textsection\ref{sec:va_on_vy} & \textsection\ref{sec:fam_on_mprs} & description \\
\hline Definition~\ref{def:va_problem}&Definition~\ref{def:fam_mprs_problem} & problem description \\
\hline Lemma~\ref{lem:va_low}&Lemma~\ref{lem:fam_low} & lowerbound  \\
\hline Lemma~\ref{lem:va_tight}&Lemma~\ref{lem:fam_tight} & properties of tight \\
\hline Proposition~\ref{prop:va_min}&Proposition~\ref{prop:fam_min} & classification minimal vineyards/mprs \\
\hline Definition~\ref{def:casedist}&Definition~\ref{def:casedist_surface} & cases A0-A4/B0-B4 \\
\hline Theorem~\ref{thm:va_opt}&Theorem~\ref{thm:optfam} & determining optimal vineyards/optimal families \\
\hline 
\end{tabular}
 } 
\end{center} The proofs of the geometric statement in this section was modeled as a  blueprint of the proof of the combinatorial in \textsection\ref{sec:va_on_vy} (we thank the anonymous referee for the notion of blueprint). The combinatorial proof served us as a guideline to a deeper understanding of the geometric one. As described in \textsection\ref{sec:fam_on_tor}, there is a translation of the vineyard problem to the family problem. Under this correspondence, the adjoint polygon (see Definition~\ref{def:adjoint_polygon})  translates into the definition of the adjoint relation for minimally polarized toric surfaces: the projective embedding defined by the interior lattice points is the embedding associated to the adjoint linear system $D+K$, where $D$ is the divisor defined by the original lattice polygon (see \cite{ful1}). So, not only the problem but also the theorem and proof translates to toric surfaces. But the so obtained theorem and proof do not use the toric structure and can be generalized to the case of arbitrary rational surfaces. 
\end{remark}

\section{Minimal degree families on rational surfaces}
 
\label{sec:fam_on_rat}
\begin{definition}
\label{def:fam_mprs_sing}
\textrm{\textbf{(optimal families and minimal degree)}}  Let $Y\subset {\textbf{P}}^r$ a rational complex surface (possibly singular) for $r\in{\textbf{Z}}_{>1}$.  Let $F\in\textrm{Fam}(Y)$. 
  The \textit{minimal rational degree}  with respect to  $Y\subset {\textbf{P}}^r$ is defined as \[v(Y)=\displaystyle\min\{~ \deg F ~~|~~ F\in\textrm{Fam} Y {\textrm{~and~}} p_g(F)=0 ~\}.\] We call $F$ an \textit{optimal family}  if and only if  $F$ is a rational family and $\deg F=v(Y)$. The \textit{set of all optimal families} on $Y\subset {\textbf{P}}^r$ is denoted by $S(Y)$. 

\end{definition}

\begin{definition}
\label{def:fam_rat_problem}
\textrm{\textbf{(rational family problem on rational surfaces)}} 
  Given a rational complex surface $Y\subset {\textbf{P}}^r$, find the minimal degree $v(Y)$ and all optimal families $S(Y)$. 
 
\end{definition}

\begin{proposition}\label{prop:optfam_rat}\textbf{\textrm{(optimal families on rational surfaces)}}  Let $Y\subset{\textbf{P}}^r$ be a rational complex surface for $r\in{\textbf{Z}}_{>1}$.  Let $\varphi_D: X\rightarrow Y$ be the minimal resolution of singularities of $Y$.  Let $D\in\textrm{Div} X$ be the pullback of hyperplane sections of $Y$.  Let $S(X,D)$ and $v(X,D)$ de defined as in Definition~\ref{def:fam_mprs}.  Let $\varphi_{D\circledast} S(X,D)=\{~ \varphi_{D\circledast}F ~~|~~ F\in S(X,D) ~\}$. 
\begin{itemize} 
\item[\textbf{a)}] We have that $(X,D)$ is a mprs. 
\item[\textbf{b)}] We have that $S(Y)=\varphi_{D\circledast} S(X,D)$ and $v(Y)=v(X,D)$. 
\end{itemize}
\begin{proof} 

\textit{Claim:} We have that \textbf{a)}. 
  It follows from $D$ being the pullback of the hyperplane sections of $Y\subset{\textbf{P}}^r$ that $D$ is nef and movable.  It follows from $\varphi_D$ being a minimal resolution that $DE>0$ for all exceptional curves $E\subset X$.  It follows from the definitions that $X$ is rational and nonsingular.

\textit{Claim:} We have that \textbf{b)}. 
  It follows from $\varphi_D$ being birational and from $\deg\varphi_D(C)=DC$ for all curves $C\subset X$ . 

 \end{proof}
\end{proposition}

\section{Examples of minimal degree families on rational surfaces}

\begin{example}\label{ex:b2}\textrm{\textbf{(case B2)}} 
 Let $F,G \in {\textbf{C}}$ be homogeneous where  
 $\deg F=d$ and $\deg G=d-1$ for $d\in{\textbf{Z}}_{\geq 4}$ and 
 $F$ and $G$ have no common multiple points. 

 Let \[ Y: F(y,z,w) - xG(y,z,w) = 0 \subset{\textbf{P}}^3, \] a complex projective surface of degree $d$.

 We consider the following birational map which parametrizes $Y$: \[ f: {\textbf{P}}^2 \dashrightarrow Y \subset {\textbf{P}}^3, \quad (s:t:u) \mapsto (F(s,t,u): sG(s,t,u):tG(s,t,u):uG(s,t,u)) \] given by polynomials of degree $d$ (also called parametric degree).

 We define $g: X \rightarrow {\textbf{P}}^2$ to be the resolution of the projective plane in the basepoints of $f$. There are $d(d-1)$ basepoints including infinitely near basepoints. We define $D\in\textrm{Cl} X$ to be associated to the resolution of $f$ which is shown in the following commutative diagram: 
\begin{center}
 \begin{tabular}
{l@{~~}l@{}l} $X$ & & \\
     $g\downarrow$ & $\searrow \varphi_{D}$ & \\
 ${\textbf{P}}^2$  & $\stackrel{f}{\dashrightarrow}$ & $Y$ 
\end{tabular}
 
\end{center}

 From $\varphi_{D}$ being a morphism it follows that $D$ is nef. Since $|D+K|$ doesn't have fixed components it follows that $(X,D)$ is a mprs (see Definition~\ref{def:mprs}). Let $E_1,\ldots,E_{d(d-1)}$ be the pullbacks of the exceptional curves resulting from blowing up the basepoints of $f$ and $L$ is the pullback of hyperplane sections of ${\textbf{P}}^2$.

 We consider the adjoint relation (see Definition~\ref{def:adjoint_relation}) \[ \ensuremath{(X,D)\stackrel{\mu}{\rightarrow}(X',D')} \] where  
 $\textrm{Cl} X={\textbf{Z}}\langle L,E_1,\ldots,E_{d(d-1)}\rangle $ and  $\textrm{Cl} X'={\textbf{Z}}\langle L\rangle $, 
 $K=-3L+E_1+\ldots+E_{d(d-1)}$ and $K'=-3L$ the canonical divisor classes, 
 $D=dL-E_1-\ldots-E_{d(d-1)}$ and $D'=(d-3)L$,  where $L^2=1, E_iE_j=-\delta_{ij}$ and $LE_i=0$ for all $i,j\in\{1,\ldots,12\}$. Let's assume that $E_1,\ldots,E_r$  correspond to the planar (not infinitly near) basepoints for $1\leq r\leq d(d-1)$.

 From $X'\cong{\textbf{P}}^2$ and Theorem~\ref{thm:optfam}.c case B2 it follows that  
 $S(X',D')=\{~F\subset L~\}$ and $v(X',D')=d-3$. 
 $S(X,D)=\{~|L-E_1|,|L-E_2|,\ldots,|L-E_r|~\}$ and $v(X,D)=d-1$. 

 For instance let $\tilde{Y}$ be an affine real representation of $Y$ where  
 $F(a_0,a_1,a_2)=a_0 (a_0 + a_2) (a_0 + 2 a_2) (a_0 + 3 a_2)$ and 
 $G(a_0,a_1,a_2)=a_1 (a_1 + a_2) (a_1 + 2 a_2)$.  The images in Figure~\ref{fig:deg4} show family members of $|L-E_1|$ on $\tilde{Y}$. From the top view in Figure~\ref{fig:deg4}.a it can be seen that the family is projected to lines through a point in the plane. The exceptional curves $E_1,\ldots,E_{12}$ are vertical lines. The family $|L-E_1|$ is given by the hyperplane sections through the vertical line corresponding to $E_1$ minus the fixed component which is the line itself. In Figure~\ref{fig:deg4}.b-d are some hyperplane sections shown corresponding to the family members.
\begin{figure}[h!]
 \begin{tabular}
{cc} {\includegraphics[width=4cm,height=4cm]{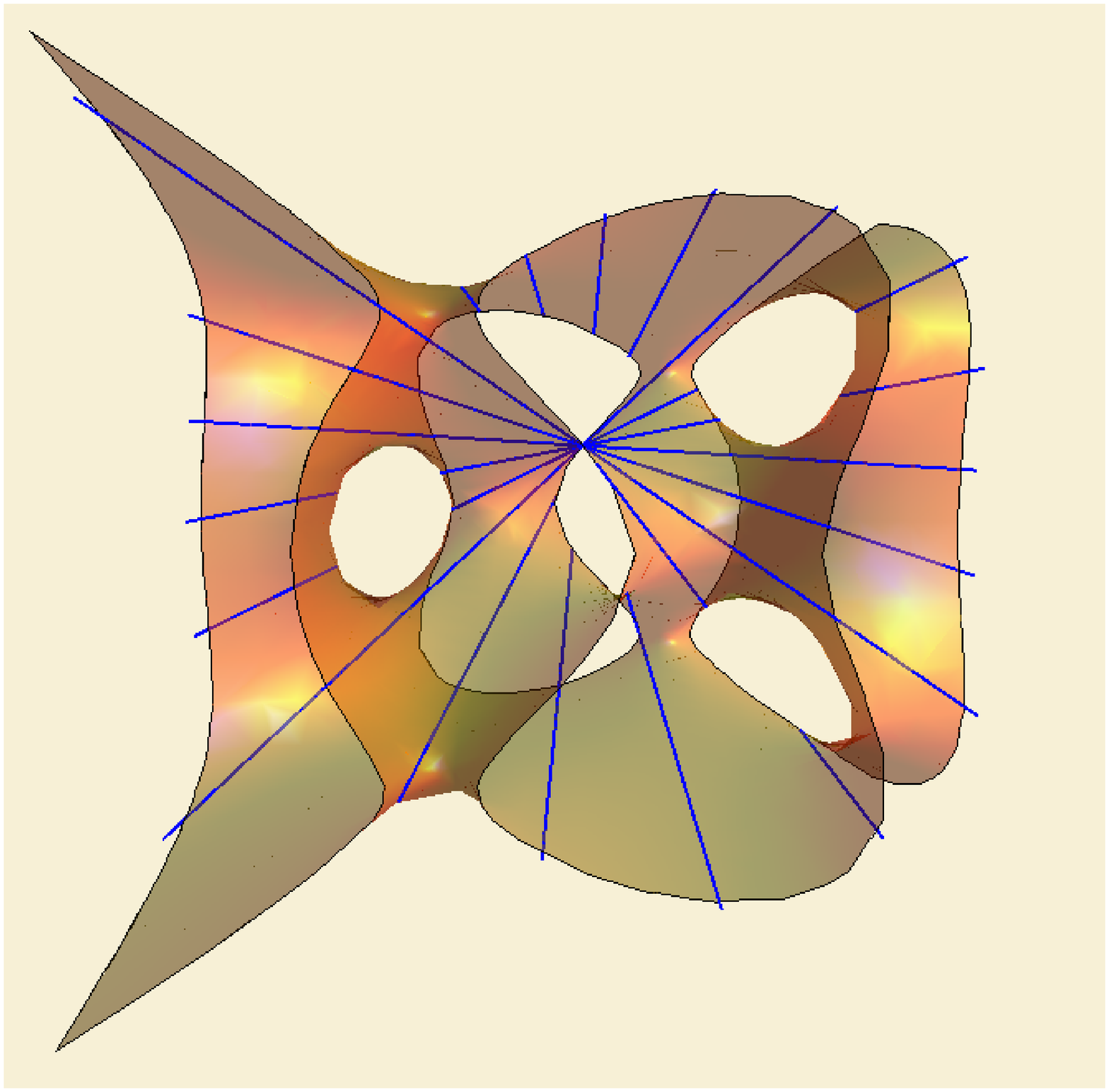}} & {\includegraphics[width=4cm,height=4cm]{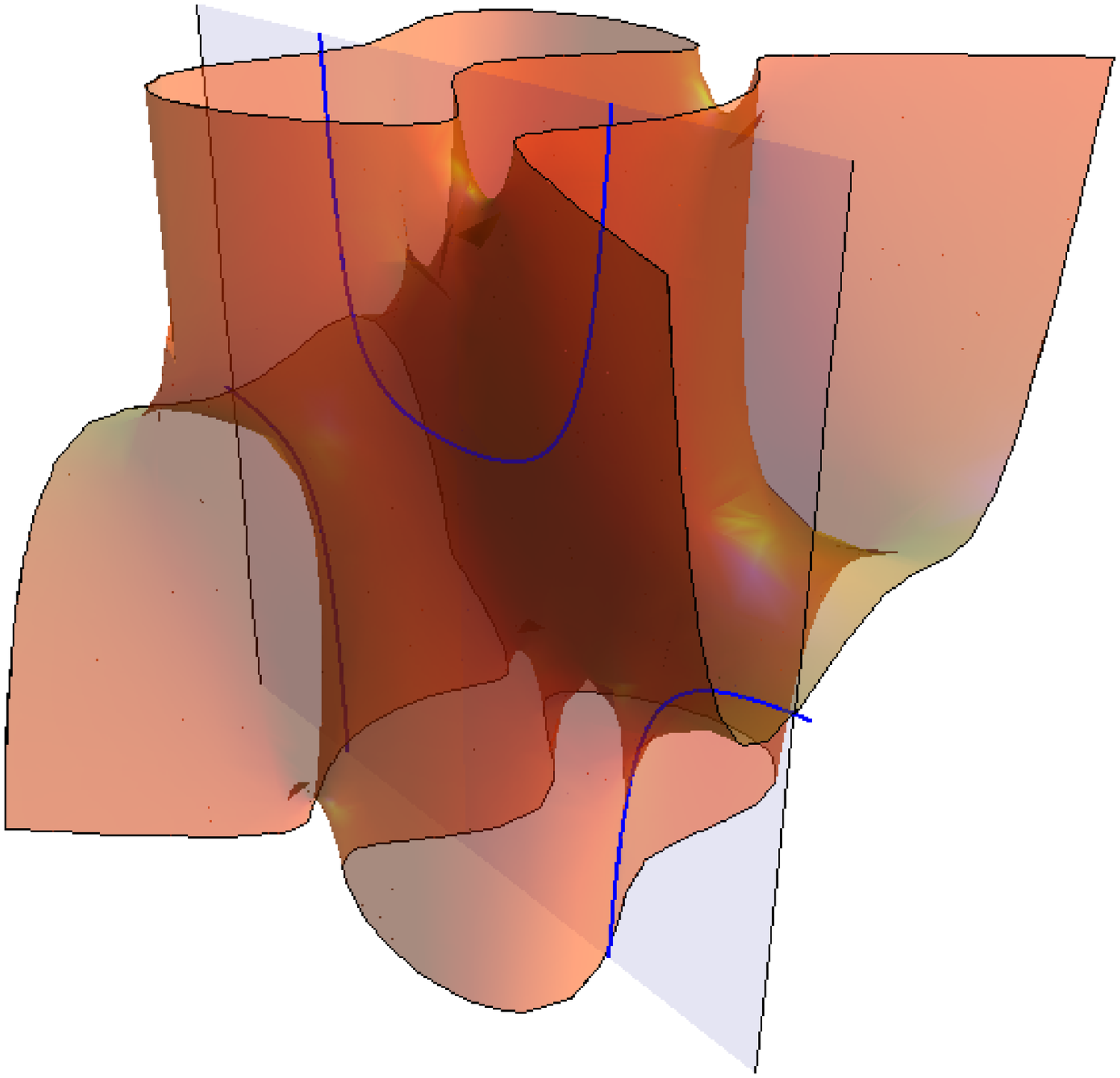}} \\
 a&b \\
 {\includegraphics[width=4cm,height=4cm]{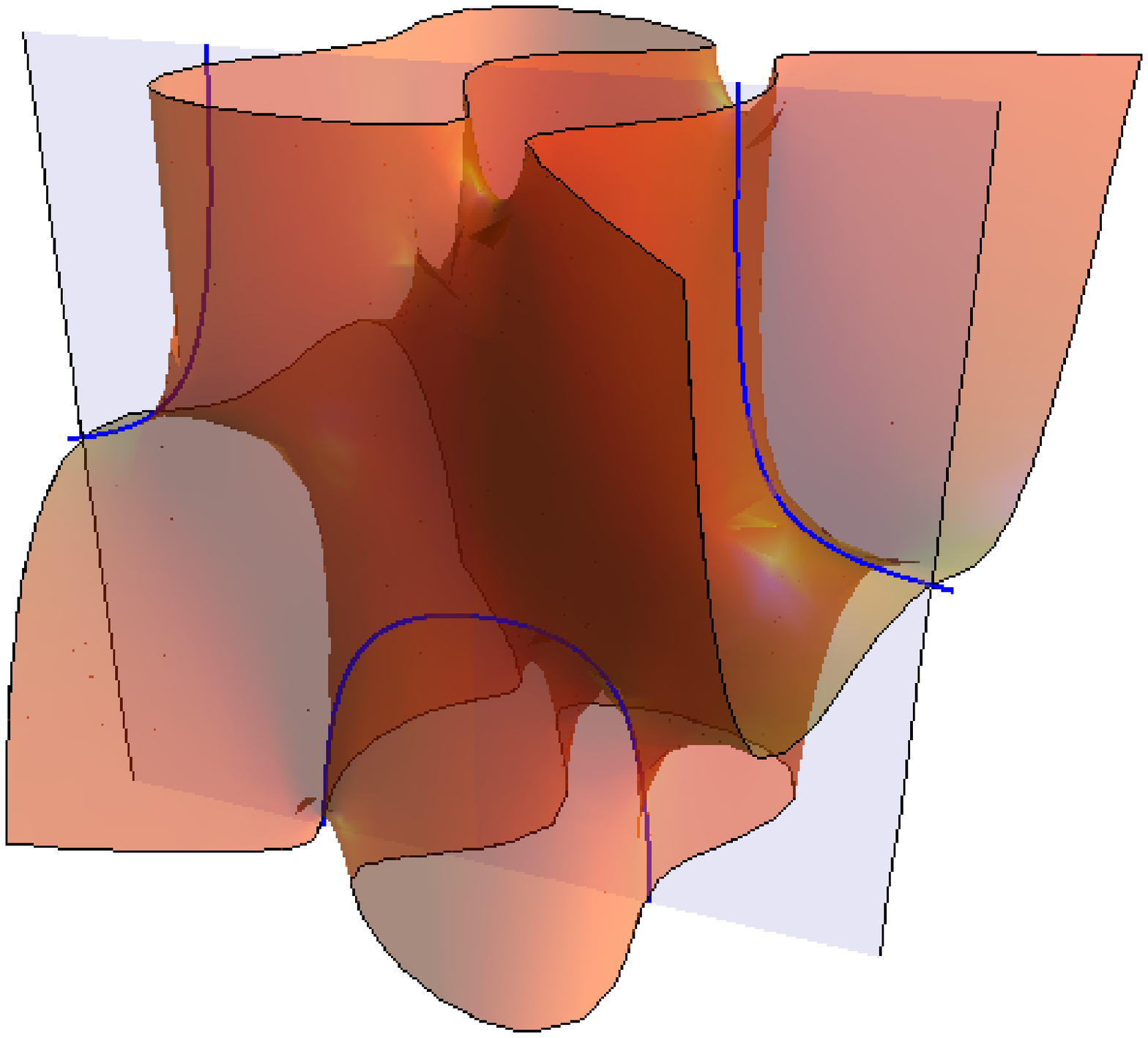}} & {\includegraphics[width=4cm,height=4cm]{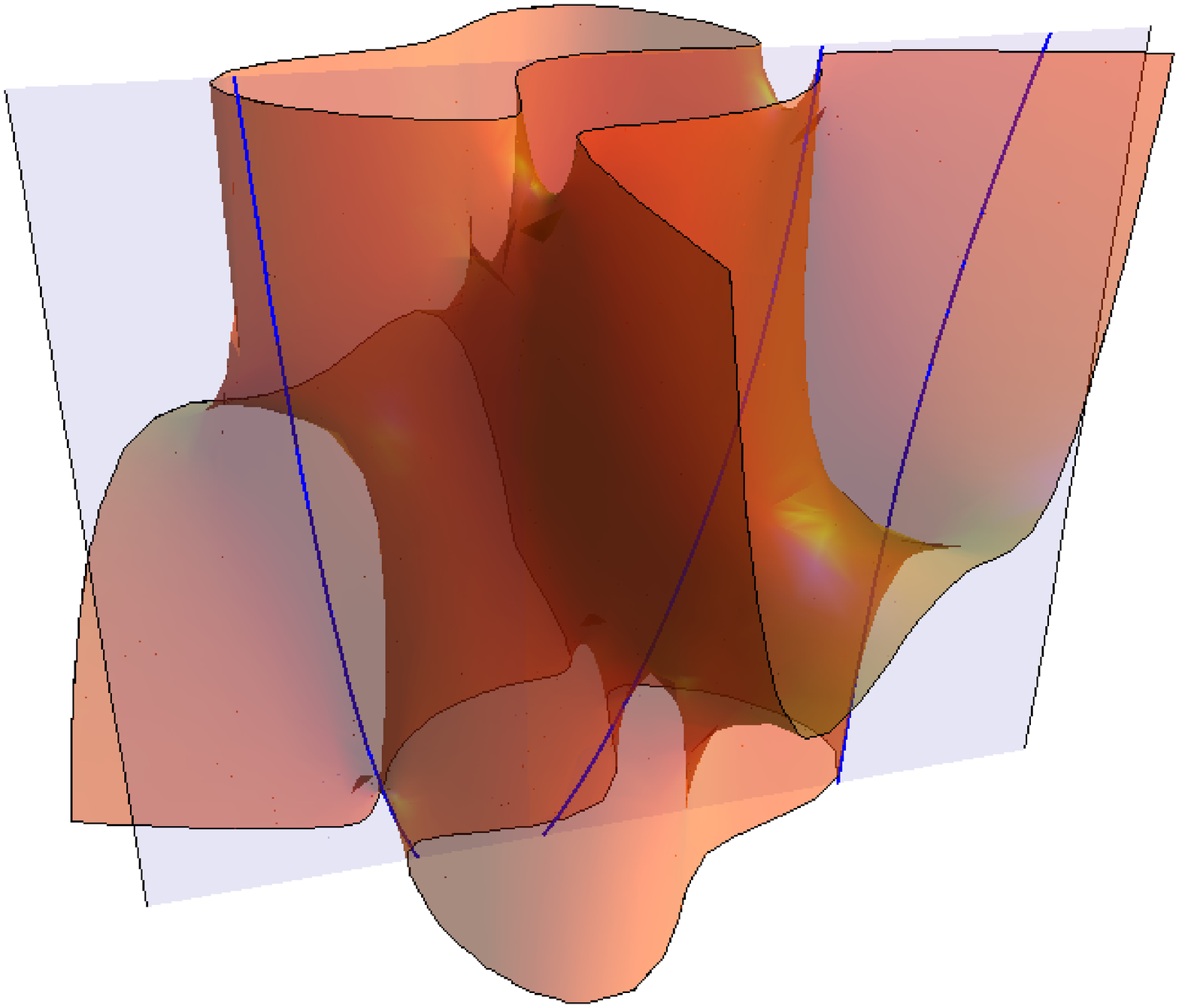}} \\
 c&d  
\end{tabular}
 \caption{Family members of an optimal family on $\tilde{Y}$ in Example~\ref{ex:b2}.} 
\label{fig:deg4} 
\end{figure}
 
\end{example}

\begin{example}\label{ex:b0_b3_b4}\textrm{\textbf{(cases B0, B3 and B4)}}

 Let $Y_0:       w^6 x^2       + 12 w^7 y       +     w^2 x^5 y    -    w^3 x^3 y^2 +   4 w^4 x^3 z     + 52 w^5 x y z   -   2 w x^4 y^2 z  +  4 w^2 x^4 z^2 +  26 w^3 x^2 y z^2 - 12 w^4 y^2 z^2 +     x^3 y^3 z^2  - 20 w x^3 y z^3 -  28 w^2 x y^2 z^3 - 24 w^3 x z^4   +  17 x^2 y^2 z^4  - 48 w x^2 z^5 - 240 w^2 y z^5     + 88 x y z^6     + 144 z^8 =0 \subset {\textbf{P}}^3 $ a complex projective surface of degree $8$.

 We consider the following birational map which parametrizes $Y_0$: $ f_0: {\textbf{P}}^2 \dashrightarrow Y_0 \subset {\textbf{P}}^3, \quad (s:t:u) \mapsto \\
(  -45 s^{15} t^{4} - 21 s^{10} t^{9} - 15 s^{5} t^{14} + t^{19} + 24 s^{18} u -    6 s^{13} t^{5} u + 20 s^{8} t^{10} u - 6 s^{3} t^{15} u - 44 s^{16} t u^{2} -    78 s^{11} t^{6} u^{2} - 60 s^{6} t^{11} u^{2} + 6 s t^{16} u^{2} -    6 s^{14} t^{2} u^{3} + 60 s^{9} t^{7} u^{3} - 30 s^{4} t^{12} u^{3} -    93 s^{12} t^{3} u^{4} - 90 s^{7} t^{8} u^{4} + 15 s^{2} t^{13} u^{4} +    60 s^{10} t^{4} u^{5} - 60 s^{5} t^{9} u^{5} - 36 s^{13} u^{6} - 60 s^{8} t^{5} u^{6} +    20 s^{3} t^{10} u^{6} + 20 s^{11} t u^{7} - 60 s^{6} t^{6} u^{7} -    15 s^{9} t^{2} u^{8} + 15 s^{4} t^{7} u^{8} - 30 s^{7} t^{3} u^{9} +    6 s^{5} t^{4} u^{10} - 6 s^{8} u^{11} + s^{6} t u^{12} \\
:  9 s^{11} t^{8} + 6 s^{6} t^{13} + s t^{18} - 12 s^{14} t^{4} u - 4 s^{9} t^{9} u +    4 s^{17} u^{2} + 18 s^{12} t^{5} u^{2} + 24 s^{7} t^{10} u^{2} + 6 s^{2} t^{15} u^{2} -    12 s^{15} t u^{3} - 12 s^{10} t^{6} u^{3} + 9 s^{13} t^{2} u^{4} +    36 s^{8} t^{7} u^{4} + 15 s^{3} t^{12} u^{4} - 12 s^{11} t^{3} u^{5} +    24 s^{9} t^{4} u^{6} + 20 s^{4} t^{9} u^{6} - 4 s^{12} u^{7} + 6 s^{10} t u^{8} +    15 s^{5} t^{6} u^{8} + 6 s^{6} t^{3} u^{10} + s^{7} u^{12} \\
:  -6 s^{18} t - 9 s^{13} t^{6} - 10 s^{8} t^{11} + s^{3} t^{16} - s^{16} t^{2} u +    10 s^{11} t^{7} u - 5 s^{6} t^{12} u - 23 s^{14} t^{3} u^{2} - 30 s^{9} t^{8} u^{2} +    5 s^{4} t^{13} u^{2} + 20 s^{12} t^{4} u^{3} - 20 s^{7} t^{9} u^{3} - 14 s^{15} u^{4} -    30 s^{10} t^{5} u^{4} + 10 s^{5} t^{10} u^{4} + 10 s^{13} t u^{5} -    30 s^{8} t^{6} u^{5} - 10 s^{11} t^{2} u^{6} + 10 s^{6} t^{7} u^{6} -    20 s^{9} t^{3} u^{7} + 5 s^{7} t^{4} u^{8} - 5 s^{10} u^{9} + s^{8} t u^{10} \\
:  -12 s^{19} - 15 s^{14} t^{5} + 10 s^{9} t^{10} + s^{4} t^{15} + 4 s^{17} t u -    20 s^{12} t^{6} u - 15 s^{15} t^{2} u^{2} + 30 s^{10} t^{7} u^{2} +    5 s^{5} t^{12} u^{2} - 40 s^{13} t^{3} u^{3} + 30 s^{11} t^{4} u^{4} +    10 s^{6} t^{9} u^{4} - 20 s^{14} u^{5} + 10 s^{12} t u^{6} + 10 s^{7} t^{6} u^{6} +    5 s^{8} t^{3} u^{8} + s^{9} u^{10} ) $ \\
 given by polynomials of degree $19$ (also called parametric degree).

 We define $g_0: X_0 \rightarrow {\textbf{P}}^2$ to be the resolution of the projective plane in the basepoints of $f_0$. There are $10$ basepoints with multiplicities $7,6,6,6,6, 6,6,6,6$ and $4$ (this example was constructed by first giving the basepoints with multiplicities and computing the implicit equation afterwards). We define $D_0\in\textrm{Cl} X_0$ to be associated to the resolution of $f_0$ which is shown in the following commutative diagram: 
\begin{center}
 \begin{tabular}
{l@{~~}l@{}l} $X_0$ & & \\
     $g_0\downarrow$ & $\searrow \varphi_{D_0}$ & \\
 ${\textbf{P}}^2$  & $\stackrel{f_0}{\dashrightarrow}$ & $Y_0$ 
\end{tabular}
 
\end{center}

 We have that $(X_0,D_0)$ is a mprs by the same argument as in Example~\ref{ex:b2}.

 We consider the adjoint chain $ \ensuremath{(X_0,D_0)\stackrel{\mu_0}{\rightarrow}(X_1,D_1)} \ensuremath{\stackrel{\mu_1}{\rightarrow}} \ensuremath{\ldots\stackrel{\mu_5}{\rightarrow}(X_6,D_6)}. $ Let $E_i$ be the pullbacks of the exceptional curves resulting from blowing up the basepoints of $f$ and $L$ is the pullback of hyperplane sections of ${\textbf{P}}^2$.

 The divisor class groups of the pairs $X_i$ for $i\in\{1,\ldots,6\}$ are generated by:  
 $\textrm{Cl} X_0=\textrm{Cl} X_1=\textrm{Cl} X_2=\textrm{Cl} X_3={\textbf{Z}}\langle L,E_1,\ldots,E_{10}\rangle $, 
 $\textrm{Cl} X_4=\textrm{Cl} X_5={\textbf{Z}}\langle L,E_1,\ldots,E_9\rangle $ and 
 $\textrm{Cl} X_6={\textbf{Z}}\langle L,E_1\rangle $,  where $L^2=1, E_iE_j=-\delta_{ij}$ and $LE_i=0$ for all $i,j\in\{1,\ldots,10\}$.

 Note that $D_{i+1}=\mu_{i*}(D_i+K_i)$, and that $\ensuremath{\textrm{dim}}|D_6+K_6|=-1$ as in Definition~\ref{def:mprs}. We can determine $D_0$ by the basepoint analysis of $f_0$. From $K_{{\textbf{P}}^2}=-3L$ (see \cite{har1}) it follows that   the canonical divisor classes $K_i \in \textrm{Cl} X_i$ are  
 $K_i=-3L+E_0+\ldots+E_{10}$ for $i\in\{0,1,2,3\}$, 
 $K_i=-3L+E_0+\ldots+E_9$ for $i\in\{4,5\}$ and 
 $K_6=-3L+E_1$. 

 We represent $D_i$ in terms of the generators of $\textrm{Cl} X_i$ for $i\in\{0,\ldots,6\}$: \\[2mm]
 {\tiny \begin{tabular}
{|l|l|l|l|l|l|l|l|l|l|l|l|} \hline       & $L$ & $-E_1$ & $-E_2$ & $-E_3$ & $-E_4$ & $-E_5$ & $-E_6$ & $-E_7$ & $-E_8$ & $-E_9$ & $-E_{10}$ \\
\hline $D_0$ & 19 & 7 & 6 & 6 & 6 & 6 & 6 & 6 & 6 & 6 & 4 \\
\hline $D_1$ & 16 & 6 & 5 & 5 & 5 & 5 & 5 & 5 & 5 & 5 & 3 \\
\hline $D_2$ & 13 & 5 & 4 & 4 & 4 & 4 & 4 & 4 & 4 & 4 & 2 \\
\hline $D_3$ & 10 & 4 & 3 & 3 & 3 & 3 & 3 & 3 & 3 & 3 & 1 \\
\hline $D_4$ &  7 & 3 & 2 & 2 & 2 & 2 & 2 & 2 & 2 & 2 & - \\
\hline $D_5$ &  4 & 2 & 1 & 1 & 1 & 1 & 1 & 1 & 1 & 1 & - \\
\hline $D_6$ &  1 & 1 & - & - & - & - & - & - & - & - & - \\
\hline 
\end{tabular}
 }

 From Theorem~\ref{thm:optfam}.a and $P=|L-E_1|$ in Proposition~\ref{prop:fam_min} it follows that  
 $S(X_6,D_6)=\{~P~\}$ and $v(X_6,D_6)=0$.  From Theorem~\ref{thm:optfam}.d case B3 and then five times B4 it follows that  
 $S(X_0,D_0)=\{~Q~\}$ and $v(X_0,D_0)=12$.  where $Q=\mu_0^\circledast\circ\ldots\circ\mu_5^\circledast P$. In terms of generators of $\textrm{Cl} X_0$ we have that $Q=|L-E_1| \in \textrm{Fam} X_0$.

 We have $g_{0\circledast}Q=F=(F_i)_{i\in I}$ with  
 $I={\textbf{P}}^1$ and $F_i=\{~ (s:\frac{i_0}{i_1}s:u) ~\}$ for $i=(i_0:i_1)\in I$.  After dehomogenization of $F$ to $s$ and $i_1$ we have  
 $I={\textbf{C}}$ and $F_i=\{~ (1:i:u) ~~|~~ u\in{\textbf{C}} ~\}^-$ for $i\in I$. 

 From Proposition~\ref{prop:optfam_rat} it follows that $S(Y_0)=\{~\varphi_{D_0\circledast}(Q)~\}$ and $v(Y_0)=v(X_0,D_0)$ where \[ \varphi_{D_0\circledast}(Q)=\varphi_{D_0\circledast}(g_0^\circledast(F))= \{~ f_0(1:i:u) ~~|~~ u\in{\textbf{C}} ~\}_{i\in I}. \] The degree of this family is $\deg_u f_0(1:i:u)=12$. Indeed this is equal to $v(X_0,D_0)=12$.

 It is remarkable that on a rational surface of degree $8$ the optimal family has degree $12$. 
 
\end{example}
 \bibliography{geometry,schicho}\paragraph{Addresses of authors:} ~\\
~\\
 Johann Radon Institute for Computational and Applied Mathematics (RICAM), Austrian Academy of Sciences,  Altenbergerstra{\ss}e 69, A-4040 Linz, Austria \\
 and \\
  Research Institute for Symbolic Computation (RISC), Johannes Kepler University, Altenbergerstrasse 69, A-4040 Linz, Austria \\
 \textbf{email:} niels.lubbes@oeaw.ac.at \\
~\\
 Johann Radon Institute for Computational and Applied Mathematics (RICAM) , Austrian Academy of Sciences , Altenbergerstra{\ss}e 69 , A-4040 Linz, Austria \\
 \textbf{email:} josef.schicho@oeaw.ac.at  
\end{document}